\documentclass[a4paper, 12pt]{article}%
\usepackage{amsmath}
\usepackage{amsfonts}
\usepackage{amssymb}
\usepackage{graphicx}%
\providecommand{\U}[1]{\protect\rule{.1in}{.1in}}
\def\inar1{INAR\,($1$)}

\begin{document}

\begin{center}
{\large Stationary underdispered \inar1 models based on the backward approach}%

\begin{tabular}
[c]{lll}%
Emad-Eldin Aly Ahmed Aly &  & Nadjib Bouzar\\
Department of Statistics and &  & Department of Mathematical Sciences,\\
Operations Research, &  & University of Indianapolis,\\
Faculty of Science, Kuwait University, &  & Indianapolis, IN 46227, USA\\
P. O. Box 5969, Safat 13060, Kuwait &  & Email: nbouzar@uindy.edu\\
Email: eealy50@gmail.com &  &
\end{tabular}

{\large Abstract}
\end{center}

Most of the stationary first-order autoregressive integer-valued (INAR(1))
models were developed for a given thinning operator using either the forward
approach or the backward approach. In the forward approach the marginal
distribution of the time series is specified and an appropriate distribution
for the innovation sequence is sought. Whereas in the backward setting, the
roles are reversed. The common distribution of the innovation sequence is
specified and the distributional properties of the marginal distribution of
the time series are studied. In this article we focus on the backward approach
in presence of the Binomial thinning operator. We establish a number of
theoretical results which we proceed to use to develop stationary INAR(1)
models with finite mean. We illustrate our results by presenting some new
\inar1 models that show underdispersion.


\newpage

\section{Introduction}

The area of integer-valued time series has attracted a lot of interest in
research and practice during the last 35 years. It started with the pioneering
work of McKenzie (1985), Al-Osh and Alzaid (1987) and McKenzie (1988). The
first models are based on the Binomial thinning operator of Steutel and van
Harn (1979). Since then, many families of new thinning-based first-order
autoregressive integer-valued models (INAR(1)) have been proposed and studied
in the literature. Al-Osh and Aly (1992), Aly and Bouzar (1994), Latour
(1998), Ristic et al. (2009) and Aly and Bouzar (2019) proposed and studied
new INAR(1) models developed by replacing the Binomial thinning operator by
other types of thinning operators. Additional references can be found in the review
article by Scotto et al. (2015).

The Binomial thinning (Steutel and van Harn (1979)) of $X,$ denoted by
$\alpha\odot X$ is defined as
\begin{equation}
\alpha\odot X=\sum_{i=1}^{X}Y_{i}, \label{bin_thin}%
\end{equation}
where $X$ is a $\mathbb{Z}_{+}$-valued random variable (rv), $\alpha\in(0,1)$
and $\{Y_{i}\}$ is a sequence of independent identically distributed (iid)
Bernoulli($\alpha$) rv's independent of $X$. The operation $\odot$
incorporates the discrete nature of the variates and acts as the analogue of
the standard multiplication used in the standard ARMA models.

Note that the Binomial thinning operator is a semigroup in the sense that%
\begin{equation}
\alpha\odot\left(  \beta\odot X\right)  =\beta\odot\left(  \alpha\odot
X\right)  =\left(  \alpha\beta\right)  \odot X. \label{semi gr}%
\end{equation}

Assume that $0<\alpha<1$, and $(\varepsilon_{t},t\geq1)$ is an iid sequence of
$\mathbb{Z}_{+}$-valued rv's. A sequence $(X_{t},t\geq0)$ of $\mathbb{Z}_{+}%
$-valued rv's is said to be an \inar1 process if
\begin{equation}
X_{t}=\alpha\odot X_{t-1}+\varepsilon_{t}  \qquad (t\geq1), \label{inar1_eq}%
\end{equation}
such that the binomial thinning $\alpha\odot X_{t-1}$ in (\ref{inar1_eq}) is
performed independently for each $t$. More precisely, we assume the existence
of an array $(Y_{i,t},\ i\geq1,\ t\geq0)$ of iid Bernoulli($\alpha$) rv's,
independent of $\{\varepsilon_{t}\}$, such that
\[
\alpha\odot X_{t-1}=\sum_{i=1}^{X_{t-1}}Y_{i,t-1}.
\]
In (\ref{inar1_eq}), $\{\varepsilon_{t}\}$ is referred to as the innovation
sequence and $\alpha$ as the coefficient of the process $\{X_{t}\}$.

If $\{X_{t}\}$ is an \inar1 process of (\ref{inar1_eq}), then the pgf
$\varphi_{X_{t}}(z)$ of $X_{t}$ and the common pgf $\Psi(z)$ of the innovation
sequence $\{\varepsilon_{t}\}$ must satisfy the functional equation%
\[
\varphi_{X_{t+1}}(z)=\varphi_{X_{t}}(1-\alpha+\alpha z)\Psi(z).
\]
If one further assumes that $\{X_{t}\}$ is stationary, then the common pgf
$\varphi_{X}(z)$ of $\{X_{t}\}$ satisfies
\begin{equation}
\varphi_{X}(z)=\varphi_{X}(1-\alpha+\alpha z)\Psi(z). \label{func_eq}%
\end{equation}
The next proposition states that equation (\ref{func_eq}) is a sufficient
condition for the existence of stationary \inar1 processes. For a proof see
for example Bouzar and Jayakumar (2008) or, in a more general setting, Aly and
Bouzar (1994).

\textbf{Proposition 1} Let $\alpha\in(0,1)$ and let $\varphi_{X}(z)$ and
$\Psi(z)$ be pgf's that satisfy the functional equation (\ref{func_eq}). Then,
there exists a stationary \inar1 process $\{X_{t}\}$ on some probability space
such that its marginal distribution and that of its innovation sequence
$\{\varepsilon_{t}\}$ have respective pgf's $\varphi_{X}(z)$ and $\Psi(z)$.

As a functional equation with two unknown pgf's ($\varphi_{X}(\cdot)$ and
$\Psi(\cdot))$, (\ref{func_eq}) can be solved in two different ways. The
forward approach: Fix a pgf, $\varphi_{X}(\cdot),$ and find $\Psi(\cdot),$ as
the solution of%
\[
\frac{\varphi_{X}(z)}{\varphi_{X}(1-\alpha+\alpha z)}=\Psi(z)
\]
provided that $\Psi(z)$ is a pgf. The backward approach: Fix a pgf,
$\Psi(\cdot),$ and find the pgf $\varphi_{X}(\cdot)$ that satisfies
(\ref{func_eq}). It can be shown that in this case%
\[
\varphi_{X}(z)=\lim_{n\longrightarrow\infty}\prod\limits_{i=0}^{n}%
\Psi(1-\alpha^{i}+\alpha^{i}z)
\]
provided that the limit exists and is a pgf.

Note that the forward approach is most useful if a researcher is interested in
a specific marginal distribution. The backward approach is most useful if a
researcher is interested in a specific distribution for the innovations
$\varepsilon_{t}$.

The forward approach has been widely used in the literature. In addition to
the above mentioned references, we cite McKenzie (2003), a review article, Joe
(1996, 2019), Zhu and Joe (2003 and 2010), and the monograph by
Wei{\ss \ (2018). }For results and references on the backward approach we note
the work of Jung et al. (2005), Pedeli and Karlis (2011), Wei\ss \ (2013),
Schweer and Wei\ss \ (2014) and Schweer and Wichelhaus (2015). We note that
Guerrero et al. (2020) proposed a third approach in which both the innovation
distribution and the marginal distribution of the stationary \inar1 process
are specified in advance. The thinning operator specific to these
distributions is identified by solving a functional equation.

In the current work, we adopt the backward approach to develop stationary
INAR(1) models with finite mean using the Binomial thinning operator. In
Section 2, we prove a number of foundational results in the context of the
backward approach. These results are then used to obtain most of the needed
distributional properties of the marginal distribution of the model under
minimal assumptions on the distribution of the innovation sequence. In
Sections 3-7, we illustrate our results of Section 2 by introducing and
studying in details the important underdispersed models when the innovations
follow the logarithmic distribution, the Bernoulli distribution, the Binomial
distribution, the Poissonian Binomial distribution and the Heine distribution,
respectively. In Section 8, we give some extensions of the previous models via convolution.

We will assume throughout the rest of this paper that $\alpha\in(0,1)$ and
that $\{f_{k}\}$ is a pmf with pgf $\Psi(z)$ such that $\Psi^{\prime
}(1)<\infty.$ We will also be using the notation $\overline{a}=1-a$ for
$a\in(0,1)$.

\section{Foundational results of the backward approach}

The main results of this Section are given in Theorems 1- 3 below.

\textbf{Theorem 1}. The function
\begin{equation}
\varphi(z)=\prod\limits_{i=0}^{\infty}\Psi(1-\alpha^{i}+\alpha^{i}z)
\label{varphi_thm1}%
\end{equation}
is a pgf. Moreover, the convergence of the infinite product is uniform over
the interval $[0,1]$ and $\varphi(z)$ satisfies
\begin{equation}
\varphi(z)=\varphi(1-\alpha+\alpha z)\Psi(z),\qquad z\in\lbrack0,1].
\label{self_dec}%
\end{equation}

\text{Proof:} First, we recall some basic results on pgf's (Feller, Vol I,
1968, is an excellent reference). Let $\{q_{k}=\sum_{i=k+1}^{\infty}f_{i}\}$
be the sequence of the tail probabilities corresponding to $\{f_{k}\}$ and
let
\[
Q(z)=\sum_{k=0}^{\infty}q_{k}z^{k},
\]
be the generating function of $\{q_{k}\}$. We have
\begin{equation}
1-\Psi(z)=(1-z)Q(z),\quad z\in\lbrack0,1] \label{tail_gf}%
\end{equation}
and
\begin{equation}
Q(1)=\sum_{k=0}^{\infty}q_{k}=\sum_{k=0}^{\infty}kf_{k}=\Psi^{\prime
}(1)<\infty. \label{finite_mean}%
\end{equation}
Define%
\[
h_{i}(z)=1-\Psi(1-\alpha^{i}+\alpha^{i}z).
\]
We have by (\ref{tail_gf}),
\[
h_{i}(z)=\alpha^{i}(1-z)Q(1-\alpha^{i}+\alpha^{i}z).
\]
Noting that $Q$ is increasing over $[0,1],0\leq1-z\leq1$, and $Q(1)$ is finite
(cf. (\ref{finite_mean})), it follows that $0\leq h_{i}(z)\leq Q(1)\alpha^{i}$
and
\[
\sum_{i=n+1}^{\infty}h_{i}(z)\leq Q(1)\sum_{i=n+1}^{\infty}\alpha^{i}.
\]
This implies $\sum_{i=n+1}^{\infty}h_{i}(z)$ converges uniformly to $0$ over
the interval $[0,1]$, which in turn implies, by Theorem 1, p. 381, in Knopp
(1990), that%
\begin{equation}
\varphi_{n+1}(z)=\prod\limits_{i=0}^{n}\Psi(1-\alpha^{i}+\alpha^{i}%
z)=\prod\limits_{i=0}^{n}(1-h_{i}(z)),n\geq0 \label{partial_prod}%
\end{equation}
converges uniformly over the interval $[0,1]$ to
\[
\varphi(z)=\prod_{i=0}^{\infty}(1-h_{i}(z))=\prod\limits_{i=0}^{\infty}%
\Psi(1-\alpha^{i}+\alpha^{i}z).
\]
Next, we show that $\lim_{z\uparrow1}\varphi(z)=1$. Define
\[
r_{n}(z)=\prod_{i=n+1}^{\infty}\Psi(1-\alpha^{i}+\alpha^{i}z).
\]
Let $\delta>0$ be arbitrary. By the uniform convergence of $\{\varphi
_{n+1}(z)\}$ to $\varphi(z)$, there exists a positive integer $N(\delta)$ such
that for any $n>N(\delta)$,
\[
\sup\limits_{z\in\lbrack0,1]}|r_{n}(z)-1|<\delta.
\]
Note that $\varphi_{n+1}(\cdot)$ of (\ref{partial_prod}) satisfies
$\varphi_{n+1}(1)=1$ and $\varphi_{n+1}(z)\leq1$. Since
\[
|\varphi(z)-1|=|\varphi_{n+1}(z)(r_{n}(z)-1)+\varphi_{n+1}(z)-1|,
\]
it follows that for any $n>N(\delta)$
\[
|\varphi(z)-1|\leq\delta+|\varphi_{n+1}(z)-1|,
\]
which in turn implies that%
\[
\limsup_{z\uparrow1}|\varphi(z)-1|=\limsup_{z\uparrow1}(1-\varphi
(z))\leq\delta+\underset{z\uparrow1}{\lim\inf}\left(  1-\varphi_{n+1}%
(z)\right)  \leq\delta.
\]
Since $\varphi(z)$ is the limit of the sequence of pgf's $\{\varphi
_{n+1}(z)\}$, we conclude by the Continuity Theorem that $\varphi(z)$ is a
pgf. Equation (\ref{self_dec}) is easily shown to hold. $\hfill\blacksquare$

Stirling numbers of the second kind, denoted by $S(r,j),$ are defined (see
Abramowitz and Stegun (1965) and Goldberg et al. (1976)) as
$S(0,0)=1,S(0,k)=S(r,0)=0$ and
\begin{equation}
S(r,j)=\frac{1}{j!}\sum_{k=0}^{j}(-1)^{j-k}\binom{j}{k}k^{r}.
\label{Stirling_2}%
\end{equation}

\textbf{Theorem 2. }We use the notation of Theorem 1. Let $\{p_{r}\}$ be the
pmf with pgf $\varphi(z)$ and let $\{f_{r}^{(i)}\}$ be the pmf with pgf
$\Psi(1-\alpha^{i}+\alpha^{i}z)$. Let $\kappa_{\lbrack r]}^{(f)}$ and
$\kappa_{\lbrack r]}^{(p)}$ be the $r$-th factorial moments of $\{f_{r}\}$ and
$\{p_{r}\}$, respectively, and assume that $\kappa_{\lbrack r]}^{(f)},r\geq1$,
are finite. Then,

\begin{enumerate}
\item $f_{r}^{(0)}=f_{r}$ and
\begin{equation}
f_{r}^{(i)}=%
\begin{cases}
f_{0}+\sum_{n=1}^{\infty}(1-\alpha^{i})^{n}f_{n}, & \mbox{if }r=0\\
\alpha^{ir}\sum_{n=r}^{\infty}\binom{n}{r}f_{n}(1-\alpha^{i})^{n-r}, &
\mbox{if }r\geq1.
\end{cases}
\label{pmf_fri}%
\end{equation}

\item
\begin{equation}
p_{r}=\lim_{k\rightarrow\infty}\Bigl(f^{(0)}\ast f^{(1)}\ast\cdots\ast
f^{(k-1)}\Bigr)_{r}, \label{pmf_varphi}%
\end{equation}
where $f^{(0)}\ast f^{(1)}\ast\cdots\ast f^{(k-1)}$ designates the $k$-factor
convolution of the pmf's $\{f_{r}^{(0)}\},\{f_{r}^{(1)}\},\cdot\cdot
\cdot,\{f_{r}^{(n-1)}\}$.

\item For every $r\geq1$, $\kappa_{\lbrack r]}^{(p)}$ and $\kappa_{r}^{(p)}$
are finite and are given by
\begin{equation}
\kappa_{\lbrack r]}^{(p)}={\frac{\kappa_{\lbrack r]}^{(f)}}{1-\alpha^{r}}}
\label{factmom_gen}%
\end{equation}
and%
\begin{equation}
\kappa_{r}^{(p)}=\sum_{j=1}^{r}S(r,j){\frac{\kappa_{\lbrack j]}^{(f)}%
}{1-\alpha^{j}}}. \label{cumfaccum_2}%
\end{equation}

\end{enumerate}

\textbf{Proof: }The proof of (\ref{pmf_fri}) is straightforward. Since
\[
\varphi(z)=\lim_{k\rightarrow\infty}\prod\limits_{i=0}^{k-1}\Psi(1-\alpha
^{i}+\alpha^{i}z),
\]
we obtain (\ref{pmf_varphi}) by the Continuity Theorem and (\ref{pmf_fri}).

Recall that the factorial cumulants, $(\kappa_{\lbrack r]},r\geq1)$, and the
cumulants, $(\kappa_{r},r\geq1)$, of a pmf are the coefficients of
${\frac{t^{r}}{r!}}$ in the power series expansions of the factorial cumulant
generating function (fcgf) $\ln\varphi(1+t)$ and the cumulant generating
function (cgf) $\ln\varphi(e^{t})$, respectively. A general formula that links
$\{\kappa_{r}\}$ and $\{\kappa_{\lbrack r]}\}$ (see Johnson et al. (2005),
Sections 1.2.7 and 1.2.8) is given by
\begin{equation}
\kappa_{r}=\sum_{j=1}^{r}S(r,j)\kappa_{\lbrack j]}, \label{cumfaccum_1}%
\end{equation}
where $S(r,j)$ are the Stirling numbers of the second kind of
(\ref{Stirling_2}).

By (\ref{varphi_thm1}),
\[
\ln\varphi(1+t)=\sum_{i=0}^{\infty}\ln\Psi(1+\alpha^{i}u)=\sum_{i=0}^{\infty
}\sum_{r=1}^{\infty}\alpha^{ir}\kappa_{\lbrack r]}^{(f)}{\frac{t^{r}}{r!}},
\]
or
\[
\ln\varphi(1+t)=\sum_{r=1}^{\infty}{\frac{\kappa_{\lbrack r]}^{(f)}}%
{1-\alpha^{r}}}{\frac{t^{r}}{r!}}.
\]
Hence we obtain (\ref{factmom_gen}). Equation (\ref{cumfaccum_2}) for the
$r$-th cumulant $\kappa_{r}^{(p)}$ follows from (\ref{factmom_gen}) and
(\ref{cumfaccum_1}). $\hfill\blacksquare$

The backward approach is based on the following result which is an immediate
consequence of (\ref{semi gr}), Theorem 1 and Proposition 1. The proof is omitted.

\textbf{Theorem 3}. Any pgf $\Psi(z)$ such that $\Psi^{\prime}(1)<\infty$,
gives rise to a stationary \inar1 process $\{X_{t}\}$ defined on some
probability space $(\Omega,\mathcal{F},P)$ and driven by equation
(\ref{inar1_eq}). Its marginal pgf is
\begin{equation}
\varphi_{X}(z)=\prod_{i=0}^{\infty}\Psi(1-\alpha^{i}+\alpha^{i}z).
\label{varphi_thm3}%
\end{equation}

\textbf{Remarks:}

\begin{enumerate}
\item Equation (\ref{factmom_gen}) was derived by Wei{\ss } (2013) using a
different approach.

\item We give a couple of possibly useful formulas for $\kappa_{\lbrack
r]}^{(p)}$ and $\kappa_{r}^{(p)}$ in terms of the sequences $(\kappa_{\lbrack
r]}^{(i)},i\geq0)$ and $(\kappa_{r}^{(i)},i\geq0)$ of the pmf's $(\{f_{r}%
^{(i)}\},i\geq0)$,
\begin{equation}
\kappa_{\lbrack r]}^{(p)}=\sum_{i=0}^{\infty}\kappa_{\lbrack r]}^{(i)}%
\quad\hbox{and}\quad\kappa_{r}^{(p)}=\sum_{i=0}^{\infty}\kappa_{r}^{(i)},
\label{alt_cumfaccum}%
\end{equation}
provided the two series converge. The proof is straightforward.

\item Provided they are finite, the first and second cumulants of a pmf, are
its mean and variance, respectively. It follows that the mean, $\mu^{(p)}$,
and the variance, $(\sigma^{(p)})^{2}$, of $\{p_{r}\}$ can be obtained from
their $\{f_{r}\}$ counterparts, $\mu^{(f)}$, $(\sigma^{(f)})^{2}$:
\begin{equation}
\mu^{(p)}={\frac{\mu^{(f)}}{1-\alpha}}\quad\hbox{and}\quad(\sigma^{(p)}%
)^{2}={\frac{(\sigma^{(f)})^{2}+\alpha\mu^{(f)}}{1-\alpha^{2}}}.
\label{mean_var}%
\end{equation}
Using (\ref{mean_var}), it is easily seen, as noted in {Wei\ss } (2013), that
$\{p_{r}\}$ of (\ref{pmf_varphi}) is underdispersed (i.e., $(\sigma^{(p)}%
)^{2}<\mu^{(p)})$ if and only if $\{f_{r}\}$ is underdispersed.

\item There are no simple formulas linking the $r$-th moment $\mu_{r}^{(p)}$
and the $r$-th factorial moment $\mu_{\lbrack r]}^{(p)}$ of $\{p_{r}\}$ to
their $\{f_{r}\}$ counterparts. However, if either $\kappa_{\lbrack r]}^{(p)}$
or $\kappa_{r}^{(p)}$ can be calculated for every $r\geq1$, then one can
compute $\mu_{r}^{(p)}$ and $\mu_{\lbrack r]}^{(p)}$ using standard formulas
that link moments and cumulants (see Johnson et al. (2005), Sections 1.2.7 and
1.2.8 and Smith (1995)).

\item The following additional results (see Al-Osh and Alzaid (1987) and
McKenzie (1988)) are needed in the sequel. An \inar1 model driven by
(\ref{inar1_eq}) is necessarily a homogeneous Markov chain with the 1-step
transition probabilities,%
\begin{equation}
P(X_{t}=k|X_{t-1}=l)=\sum_{j=0}^{\min(l,k)}\binom{l}{j}\alpha^{j}%
(1-\alpha)^{l-j}P(\varepsilon=k-j). \label{MarCh}%
\end{equation}
The $k$-step-ahead version of (\ref{inar1_eq}) for $k\geq1$ is given by
\begin{equation}
X_{t+k}\overset{D}{=}\alpha^{k}\circ X_{t}+\sum_{j=1}^{k}\alpha^{j-1}%
\circ\varepsilon_{t+k-j+1}. \label{inar1_eq2}%
\end{equation}
It follows from (\ref{inar1_eq2}) that the conditional pgf of $X_{t+k}$ given
$X_{t}$ satisfies%
\begin{equation}
\varphi_{X_{t+k}|X_{t}}(z)=\left(  1-\alpha^{k}+\alpha^{k}z\right)  ^{X_{t}%
}\times\prod\nolimits_{i=0}^{k-1}\Psi(1-\alpha^{i}+\alpha^{i}z).
\label{k_step_pgf}%
\end{equation}

\end{enumerate}

\section{Stationary \inar1 models with logarithmic innovations}

We start out by recalling a few facts about the logarithmic distribution (see
Johnson et al., 2005). The pmf of the logarithmic($p$) distribution is given
by $f_{r}={\frac{p^{r}}{-r\ln\overline{p}}},r\geq1,$ where $p\in(0,1)$. Its
pgf, mean, variance and dispersion index are given respectively by%
\begin{equation}
\Psi(z)={\frac{\ln(1-pz)}{\ln\overline{p}},} \label{logpgf}%
\end{equation}%
\[
\mu^{(f)}=-{\frac{p}{\overline{p}\ln\overline{p}},}\quad(\sigma^{(f)}%
)^{2}=-{\frac{p+\ln\overline{p}}{(\overline{p}\ln\overline{p})^{2}}}%
\quad\text{and}\quad I^{(f)}={\frac{p+\ln\overline{p}}{\overline{p}%
\ln\overline{p}}}.
\]
Note that the logarithmic distribution is underdispersed if $p<1-1/e,$
equidispersed if $p=1-1/e$ and overdispersed if $p>1-1/e$.

The factorial moments of $\{f_{r}\}$ are
\begin{equation}
\mu_{[r]}^{(f)}=-{\frac{p^{r}(r-1)!}{(1-p)^{r}\ln\overline{p}}} \qquad(r\ge1).
\label{factmom_log}%
\end{equation}

Recall that the moments $\{\mu_{r}\}$ of a random variable can be derived from
their factorial counterparts $\{\mu_{\lbrack r]}\}$ via the equation (see
Sections 1.2.7 and 1.2.8 in Johnson et al. (2005))
\begin{equation}
\mu_{r}=\sum_{j=1}^{r}S(r,j)\mu_{\lbrack j]} \qquad(r\ge1), \label{momfacmom}%
\end{equation}
with $S(r,j)$ of (\ref{Stirling_2}). By, (\ref{factmom_log}) and
(\ref{momfacmom}), the moments $\left(  \mu_{r}^{(f)}\right)  $ of a
logarithmic distribution are given by
\begin{equation}
\mu_{r}^{(f)}=-{\frac{1}{\ln\overline{p}}}\sum_{j=1}^{r}S(r,j){\frac
{p^{j}(j-1)!}{(1-p)^{j}}} \qquad(r\ge1). \label{mom_log}%
\end{equation}

\textbf{Lemma 1.} Let $p\in(0,1)$ and let $\{f_{r}\}$ be the pmf of a
logarithmic($p$) distribution with pgf $\Psi(z)$ of (\ref{logpgf}). Then for
every $i\geq0$, the pmf, $\{f_{r}^{(i)}\}$ of (\ref{pmf_fri}) (with pgf
$\Psi(1-\alpha^{i}+\alpha^{i}z)$) is a two-mixture of the Dirac measure
$\delta_{0}$ sitting at $0$ and the logarithmic($q_{i}$) distribution, with
$q_{i}={\frac{p\alpha^{i}}{1-p(1-\alpha^{i})}}$, and with respective mixing
probabilities $b_{i}=1-{\frac{\ln\overline{q}_{i}}{\ln\overline{p}}}$ and
$1-b_{i}={\frac{\ln\overline{q}_{i}}{\ln\overline{p}}}$, i.e.,
\begin{equation}
f_{r}^{(i)}=b_{i}\delta_{0}(\{r\})+(1-b_{i}){\frac{q_{i}^{r}}{-r\ln
\overline{q}_{i}}}. \label{pmf_fri_log}%
\end{equation}
Note that, $q_{0}=p,b_{0}=0$ and $f_{r}^{(0)}=f_{r}$. Moreover, the $k$-factor
convolution of the pmf's $\{f_{r}^{(0)}\},\{f_{r}^{(1)}\},\cdots
,\{f_{r}^{(k-1)}\},k\geq2$, is a finite mixture of convolutions of logarithmic
distributions, namely,
\begin{equation}
\Bigl(f^{(0)}\ast f^{(1)}\ast\cdots\ast f^{(k-1)}\Bigr)_{r}=C_{I,0}%
\,g_{r}^{(0)}+\sum_{l=1}^{k-1}\sum\limits_{\mathbf{j}\in\mathbf{J}_{l}%
}C_{\mathbf{j},l}\Bigl(g^{(0)}\ast g^{(j_{1})}\ast g^{(j_{2})}\cdots
g^{(j_{l})}\Bigr)_{r}, \label{pmf_finite_mix_log}%
\end{equation}
where $\{g_{r}^{(j)}\}$ is the pmf of the logarithmic($q_{j}$),
$I=\{1,2,\cdots,k-1\}$, $\mathbf{J}_{l}$ is the collection of ordered
$l$-tuples $\mathbf{j}=(j_{1},j_{2},\cdots,j_{l})$, $1\leq j_{1}<j_{2}%
<\cdots<j_{l}\leq k-1$ and $\mathbf{j}_{u}=\{j_{1},j_{2},\cdots,j_{l}\}$ is
the corresponding unordered $l$-tuple. The mixing probabilities are
\begin{equation}
C_{I,0}=\prod\limits_{j=1}^{k-1}b_{j}\quad\hbox{and}\quad C_{\mathbf{j}%
,l}=\Bigl(\prod_{j\in I\setminus\mathbf{j}_{u}}b_{j}\Bigr)\Bigl(\prod
_{h=1}^{l}(1-b_{j_{h}})\Bigr). \label{mix_prob_log}%
\end{equation}

\textbf{Proof:} If $i=0$, (\ref{pmf_fri_log}) is true since $\{f_{r}%
^{(0)}\}=\{f_{r}\}$. Assume $i\geq1$. By (\ref{pmf_fri}),
\[
f_{0}^{(i)}={\frac{-1}{\ln\overline{p}}}\sum_{n=1}^{\infty}{\frac{\left(
p(1-\alpha^{i})\right)  }{n}}={\frac{\ln(1-p(1-\alpha^{i}))}{\ln\overline{p}}%
},
\]
and for $r\geq1$,
\[
f_{r}^{(i)}=-{\frac{(p\alpha^{i})^{r}}{r\ln\overline{p}}}\sum_{n=r}^{\infty
}\binom{n-1}{r-1}\left(  p(1-\alpha^{i})\right)  ^{n-r}.
\]
Using the power series expansion
\begin{equation}
(1-t)^{-r-1}=\sum_{n=r}^{\infty}\binom{n}{r}t^{n-r}, \label{pow_series}%
\end{equation}
with $t=p(1-\alpha^{i})$, it follows that
\[
f_{r}^{(i)}=-{\frac{1}{r\ln\overline{p}}}\Biggl[{\frac{p\alpha^{i}%
}{1-p(1-\alpha^{i})}}\Biggr]^{r} \qquad(r\geq1).
\]
Setting $q_{i}={\frac{p\alpha^{i}}{1-p(1-\alpha^{i})}}$, it is easily verified
that $f_{r}^{(i)}$ satisfies (\ref{pmf_fri_log}). The second part of the Lemma
and equations (\ref{pmf_finite_mix_log}) and (\ref{mix_prob_log}) are proved
by a tedious but straightforward induction argument. The details are omitted.
$\hfill\blacksquare$

\textbf{Theorem 4.} Let $\{X_{t}\}$ be the stationary \inar1 process driven by
(\ref{inar1_eq}) and with a logarithmic($p$) innovation sequence for some
$p\in(0,1)$. Then,

\begin{enumerate}
\item[(i)] the marginal distribution $\{p_{r}\}$ of $\{X_{t}\}$ is given by
(\ref{pmf_varphi}), where $f^{(0)}\ast f^{(1)}\ast\cdots\ast f^{(k-1)}$ is
described by equations (\ref{pmf_finite_mix_log}) and (\ref{mix_prob_log}).

\item[(ii)] the marginal pgf $\varphi(z)$ of $\{X_{t}\}$ admits the
representation
\begin{equation}
\varphi(z)=\prod_{i=0}^{\infty}\Bigl[1-{\frac{1}{\ln\overline{p}}}\cdot
\ \ln{\frac{\overline{q}_{i}}{1-q_{i}z}}\Bigr]. \label{pgf_inar1_log}%
\end{equation}

\end{enumerate}

\textbf{Proof: }Part (i) is a direct consequence of (\ref{pmf_varphi}) and
Lemma 1. For (ii), we note the pgf of $\{f_{r}^{(i)}\}$ of (\ref{pmf_fri_log})
is
\[
\Psi(1-\alpha^{i}+\alpha^{i}z)=1-{\frac{1}{\ln\overline{p}}}\cdot\ \ln
{\frac{\overline{q}_{i}}{1-q_{i}z}},
\]
which implies that the pgf of $f^{(0)}\ast f^{(1)}\ast\cdots\ast f^{(k-1)}$
is
\[
\varphi_{k}(z)=\prod_{i=0}^{k-1}\Bigl[1-{\frac{1}{\ln\overline{p}}}\cdot
\ \ln{\frac{\overline{q}_{i}}{1-q_{i}z}}\Bigr],
\]
which in turn implies (\ref{pgf_inar1_log}).$\hfill\blacksquare$

We provide additional properties for a stationary \inar1 process $\{X_{t}\}$
with a logarithmic($p$) innovation sequence.

By (\ref{MarCh}), the 1-step transition probability is given by
\[
P(X_{t}=k|X_{t-1}=l)=-{\frac{p^{k}}{\ln\overline{p}}}\sum_{j=0}^{\min
(l,k-1)}\binom{l}{j}{\frac{(\alpha/p)^{j}\,\overline{\alpha}^{l-j}}{k-j}}.
\]

By (\ref{k_step_pgf}), the conditional pgf of $X_{t+k}$ given $X_{t}$
satisfies%
\[
\varphi_{X_{t+k}|X_{t}}(z)=\left(  1-\alpha^{k}+\alpha^{k}z\right)  ^{X_{t}%
}\times\prod_{i=0}^{k-1}\Bigl[1-{\frac{1}{\ln\overline{p}}}\cdot\ \ln
{\frac{\overline{q}_{i}}{1-q_{i}z}}\Bigr].
\]
Therefore, given $X_{t}=n$, the distribution of $X_{t+k}$ is the convolution
of a Binomial($n,\alpha^{k}$) distribution and the finite mixture of
convolutions of logarithmic distributions described by
(\ref{pmf_finite_mix_log}) and (\ref{mix_prob_log}).

By (\ref{mean_var}) and (\ref{factmom_log}), the mean ($\mu_{X}$) and variance
($\sigma_{X}^{2}$) of the marginal distribution of $\{X_{t}\}$ are given by
\[
\mu_{X}=-{\frac{p}{\overline{p}(1-\alpha)\ln\overline{p}}}\quad\hbox{and}\quad
\sigma_{X}^{2}=-{\frac{p}{\overline{p}(1-\alpha^{2})\ln\overline{p}}%
}\Bigl(\alpha+{\frac{(p+\ln\overline{p})}{\overline{p}\ln\overline{p}}%
}\Bigr).
\]
Note that the distribution of $X_{t}$ is underdispersed if $p<1-1/e.$

The moments and cumulants of the marginal distribution of $\{X_{t}\}$ are
computed as follows:

\begin{enumerate}
\item Compute the $r$-th cumulant $\kappa_{r}^{(f)}$ of $\varepsilon_{t}$
using the formula (due to Smith, 1995)
\begin{equation}
\kappa_{r}^{(f)}=\mu_{r}^{(f)}-\sum_{i=1}^{r-1}\binom{r-1}{i}\kappa
_{r-i}^{(f)}\mu_{i}^{(f)}, \label{cum_from_mom}%
\end{equation}
along with (\ref{mom_log}) (recall $\kappa_{1}^{(f)}=\mu_{\varepsilon}$ and
$\kappa_{2}^{(f)}=\sigma_{\varepsilon}^{2}$).

\item Compute the $r$-th factorial cumulant $\kappa_{\lbrack r]}^{(f)}$ of
$\varepsilon_{t}$ using the formula (see Johnson et al. (2005), Sections 1.27
and 1.2.8)
\begin{equation}
\kappa_{\lbrack r]}^{(f)}=\sum_{j=0}^{r}s(r,j)\kappa_{j}^{(f)},
\label{factcum_cum}%
\end{equation}
where $s(r,j)$ is the Stirling number of the first kind satisfying the
recurrence relation
\begin{equation}
s(r+1,j)=s(r,j-1)-ns(r,j), \label{Stirling_1}%
\end{equation}
with $s(n,0)=0$ and $s(1,1)=1$ (see Johnson et al., 2005).

\item Compute the $r$-th factorial cumulant $\kappa_{\lbrack r]}^{(p)}$ of
$X_{t}$ using (\ref{factmom_gen}). Use the approach described in Theorem 2 to
obtain the cumulant, moment, and factorial moment of $X_{t}$.
\end{enumerate}

\section{Stationary \inar1 models with\textbf{ } Bernoulli innovations}

We start out by recalling a result and a definition that will be needed in the
sequel. Let $q,c\in(0,1)$ and $m\geq2$. Kemp (1987) (see also Johnson et al.,
2005, p. 467) introduced and studied the Poissonian Binomial ($m, q, c$)
distribution as the distribution of a finite convolution of Bernoulli($cq^{i}%
$) distributions, $i=0,1,2,\cdots,m-1$ with pgf
\begin{equation}
\Psi(z)=\prod_{i=0}^{m-1}(1-c q^{i}(1-z)) \label{pgf_Po_Bin}%
\end{equation}
and pmf
\begin{equation}
q_{r}(m,q,c)=\sum_{k=r}^{m}(-1)^{k-r}\binom{k}{r}c^{k}q^{\binom{k}{2}}%
\prod_{l=0}^{k-1}{\frac{1-q^{m-l}}{1-q^{l+1}}}\quad,r=0,1,\cdots,m.
\label{pmf_Po_Bin}%
\end{equation}

A distribution on $\mathbb{Z}_{+}$ is said to have a discrete pseudo compound
Poisson distribution, $PCPD(\lambda,\{a_{k}\}),$ if its pgf can be written as
\begin{equation}
P(z)=\exp\left\{  \lambda\left(  \sum_{k=1}^{\infty}a_{k}\left(
z^{k}-1\right)  \right)  \right\}  \label{pcpd_pgf}%
\end{equation}
for some $\lambda>0$ and some sequence of real numbers $(a_{k},k\geq1)$ such
that $\sum\limits_{k=1}^{\infty}a_{k}=1$ and $\sum\limits_{k=1}^{\infty}%
|a_{k}|<\infty$.

\textbf{Theorem 5. }Let $\{X_{t}\}$ be the stationary \inar1 process driven by
(\ref{inar1_eq}) and with a Bernoulli($p$) innovation sequence for some
$p\in(0,1)$. Then,

\begin{enumerate}
\item the marginal pmf $\{p_{r}\}$ of $\{X_{t}\}$ is the weak limit of
Poissonian Binomial($n,\alpha,p$) (see (\ref{pgf_Po_Bin}) and
(\ref{pmf_Po_Bin})) as $n\rightarrow\infty$ and is given by
\begin{equation}
p_{r}=\lim_{n\rightarrow\infty}q_{r}(n,\alpha,p)=\sum_{k=r}^{\infty}%
(-1)^{k-r}\binom{k}{r}{\frac{p^{k}\alpha^{\binom{k}{2}}}{\prod_{l=1}%
^{k}(1-\alpha^{l})}},r\geq0. \label{pmf_Bern}%
\end{equation}

\item the tail probabilities $P(X_{t}\geq r)=\sum\limits_{j=r}^{\infty}p_{j}$
of $X_{t}$ are obtained by the formula
\begin{equation}
P(X_{t}\geq r)=\sum_{k=r}^{\infty}(-1)^{k-r}\binom{k-1}{r-1}{\frac{p^{k}%
\alpha^{\binom{k}{2}}}{\prod_{l=1}^{k}(1-\alpha^{l})},}\quad r\geq1.
\label{tail_Bern}%
\end{equation}

\item the marginal pgf $\varphi_{X}(z)$ of $\{X_{t}\}$ admits two useful
representations:
\begin{equation}
\varphi_{X}(z)=1+\sum_{n=1}^{\infty}{\frac{p^{n}(z-1)^{n}\alpha^{\binom{n}{2}%
}}{\prod_{l=1}^{n}(1-\alpha^{l})}} \label{Bern1}%
\end{equation}
and
\begin{equation}
\varphi_{X}(z)=\exp\left\{  -\sum_{n=1}^{\infty}{\frac{p^{n}}{n(1-\alpha^{n}%
)}}(1-z)^{n}\right\}  . \label{Bern2}%
\end{equation}

\item If $p<1/2$, then $\{X_{t}\}$ has a $PCPD(\lambda,\{a_{k}(p,\alpha)\})$
marginal distribution, where
\begin{equation}
\lambda=\sum_{n=1}^{\infty}{\frac{p^{n}}{n(1-\alpha^{n})}}\quad\hbox{and}\quad
a_{n}(p,\alpha)=\frac{\left(  -1\right)  ^{n+1}}{\lambda}\sum_{j=n}^{\infty
}\binom{j}{n}{\frac{p^{j}}{j(1-\alpha^{j})}.} \label{pcpd_2}%
\end{equation}

\end{enumerate}

\textbf{Proof:} It is long. We defer it to the Appendix Section.

Next, we establish several properties of the marginal distribution of the
stationary \inar1 process $\{X_{t}\}$ with Bernoulli($p$) innovations.

By (\ref{MarCh}), the 1-step transition probability is given by%

\begin{equation}
P(X_{t}=k|X_{t-1}=l)=\left\{
\begin{array}
[c]{cc}%
0, & k>l+1\\
p\alpha^{k-1}, & k=l+1\\
\alpha^{k-1}\overline{\alpha}^{l-k}\left\{  p\binom{l}{k-1}\overline{\alpha
}+\overline{p}\binom{l}{k}\alpha\right\}  , & k\leq l
\end{array}
\right.  . \label{tranber}%
\end{equation}

By (\ref{k_step_pgf}), the conditional pgf of $X_{t+k}$ given $X_{t}$
satisfies%
\[
\varphi_{X_{t+k}|X_{t}}(z)=\left(  1-\alpha^{k}+\alpha^{k}z\right)  ^{X_{t}%
}\times\prod_{i=0}^{k-1}(1-p\alpha^{i}(1-z)).
\]
Therefore, given $X_{t}=n$, the distribution of $X_{t+k}$ is the convolution
of a Binomial($n,\alpha^{k}$) distribution and the Poissonian Binomial
($k,\alpha,p$) distribution of (\ref{pmf_Po_Bin}).

By (\ref{mean_var}), the mean, the variance and the index of dispersion of
$X_{t}$ are
\[
\mu_{X}=\frac{p}{1-\alpha},\quad\sigma_{X}^{2}={\frac{p(1-p)+\alpha
p}{1-\alpha^{2}}\quad\hbox{and}}\quad I_{X}=1-{\frac{p}{1+\alpha}.}%
\]
Clearly, $\{X_{t}\}$ is underdispersed.

We derive the factorial moments $(\mu_{\lbrack r]},n\geq1)$ of $X_{t}$. Using
the version (\ref{Bern1}) of $\varphi_{X}(z)$, we deduce that
\[
\varphi_{X}(1+t)=1+\sum_{r=1}^{\infty}{\frac{r!p^{r}\alpha^{\binom{r}{2}}%
}{\prod_{i=1}^{r}(1-\alpha^{i})}}\cdot{\frac{t^{r}}{r!}}.
\]
Since the series converges everywhere, the factorial moments of $X_{t}$ of all
orders are finite and are given by
\begin{equation}
\mu_{\lbrack r]}=\frac{r!p^{r}\alpha^{\binom{r}{2}}}{\prod\limits_{i=1}%
^{r}(1-\alpha^{i})} \qquad(r\geq1). \label{fact_mom_Bern}%
\end{equation}
For example, we have
\begin{equation}
\mu_{\lbrack1]}={\frac{p}{1-\alpha}}\quad\text{and}\quad\mu_{\lbrack2]}%
={\frac{2\alpha p^{2}}{(1-\alpha)(1-\alpha^{2})}}. \label{facmom}%
\end{equation}

The moments $\mu_{r}$ of $X_{t}$, $r\geq1$, can be obtained from their
factorial counterparts via the formula (see (\ref{momfacmom}) and
(\ref{Stirling_2}))
\begin{equation}
\mu_{r}=\sum_{j=1}^{r}S(r,j)\frac{j!p^{j}\alpha^{\binom{j}{2}}}{\prod
\limits_{i=1}^{j}(1-\alpha^{i})} \qquad(r\ge1). \label{mom_Bern}%
\end{equation}

By (\ref{Bern2}), the fcgf of $X_{t}$ is given by
\[
\ln\varphi_{X}(1+t)=\sum_{r=1}^{\infty}{\frac{(-1)^{r+1}(r-1)!p^{r}}%
{(1-\alpha^{r})}}\cdot{\frac{t^{r}}{r!}}.
\]
Since the series above converges everywhere, the factorial cumulants of all
orders are finite and are given by
\begin{equation}
\kappa_{\lbrack r]}=(-1)^{r+1}{\frac{(r-1)!p^{r}}{(1-\alpha^{r})}}\qquad
(r\geq1). \label{fact_cum_Bern}%
\end{equation}

The cumulants of $X_{t}$ $(\kappa_{r},r\geq1)$ can be obtained from the
factorial cumulants via the formula (\ref{cumfaccum_1}) and (\ref{Stirling_2}%
)
\begin{equation}
\kappa_{r}=\sum_{j=0}^{r}S(r,j)(-1)^{j+1}{\frac{(j-1)!p^{j}}{(1-\alpha^{j})}}
\qquad(r\ge1). \label{cumfaccum_Bern}%
\end{equation}

For example,
\[
\kappa_{1}={\frac{p}{1-\alpha},}%
\]%
\[
\kappa_{2}=-{\frac{p^{2}}{1-\alpha^{2}}}+{\frac{p}{1-\alpha},}%
\]%
\[
\kappa_{3}=2{\frac{p^{3}}{1-\alpha^{3}}}-3{\frac{p^{2}}{1-\alpha^{2}}}%
+{\frac{p}{1-\alpha}}%
\]
and%
\[
\kappa_{4}=-6{\frac{p^{4}}{1-\alpha^{4}}}+12{\frac{p^{3}}{1-\alpha^{3}}%
}-7{\frac{p^{2}}{1-\alpha^{2}}}+{\frac{p}{1-\alpha}.}%
\]

\section{Stationary \inar1 models with\textbf{ } Binomial\textbf{
}innovations}

The treatment is essentially similar to the Bernoulli case ($m=1$). We
summarize the main results with minimal justifications for the most part.

\textbf{Theorem 6. }Let $\{X_{t}\}$ be the stationary \inar1 process driven by
(\ref{inar1_eq}) and with a Binomial($m,p$) innovation sequence for some
positive integer $m$ and some $p\in(0,1)$. Then

\begin{enumerate}
\item the marginal pmf $\{p_{r}\}$ of $\{X_{t}\}$ is the $m$-fold convolution
of the marginal distribution (\ref{pmf_Bern}) of the \inar1 process with a
Bernoulli($p$) innovation, or
\begin{equation}
p_{r}=\Biggl[\sum_{k=r}^{\infty}(-1)^{k-r}\binom{k}{r}{\frac{p^{k}%
\alpha^{\binom{k}{2}}}{\prod_{l=1}^{k}(1-\alpha^{l})}}\Biggr]^{\ast m}%
\quad(r\geq0). \label{pmf_Bin}%
\end{equation}

\item the marginal pgf $\varphi_{X}(z)$ of $\{X_{t}\}$ admits two
representations:
\begin{equation}
\varphi_{X}(z)=\Biggl[1+\sum_{n=1}^{\infty}{\frac{p^{n}(z-1)^{n}\alpha
^{\binom{n}{2}}}{\prod_{l=1}^{n}(1-\alpha^{l})}}\Biggr]^{m} \label{Bin1}%
\end{equation}

and%

\begin{equation}
\varphi_{X}(z)=\exp\Big\{-m\sum_{n=1}^{\infty}{\frac{p^{n}}{n(1-\alpha^{n})}%
}(1-z)^{n}\Bigr\}. \label{Bin2}%
\end{equation}

\item If $p<1/2$, then $\{X_{t}\}$ has a $PCPD(\lambda,\{a_{k}(p,\alpha)\})$
marginal distribution, where
\begin{equation}
\lambda=m\sum_{n=1}^{\infty}{\frac{p^{n}}{n(1-\alpha^{n})}}\quad
\hbox{and}\quad a_{n}(p,\alpha)=\frac{\left(  -1\right)  ^{n+1}}{\lambda}%
\sum_{j=n}^{\infty}\binom{j}{n}{\frac{p^{j}}{j(1-\alpha^{j})}.} \label{pcpd_3}%
\end{equation}

\end{enumerate}

\textbf{Proof:} Omitted.

We proceed to give additional properties for the marginal distribution of the
stationary \inar1 process $\{X_{t}\}$ with Binomial$(m,p)$ innovations.

By (\ref{MarCh}), the 1-step transition probability is given by%

\begin{equation}
P(X_{t}=k|X_{t-1}=l)=p^{k}\overline{p}^{m-k}\overline{\alpha}^{l}\sum
_{j=\max(k-m,0)}^{\min(l,k)}\binom{l}{j}\binom{m}{k-j}\Bigl({\frac
{\alpha\overline{p}}{p\overline{\alpha}}}\Bigr)^{j}. \label{TranBin}%
\end{equation}

By (\ref{k_step_pgf}), the conditional pgf of $X_{t+k}$ given $X_{t}$
satisfies%
\[
\varphi_{X_{t+k}|X_{t}}(z)=\left(  1-\alpha^{k}+\alpha^{k}z\right)  ^{X_{t}%
}\times\Bigl[\prod_{i=0}^{k-1}(1-p\alpha^{i}(1-z))\Bigr]^{m}.
\]
Therefore, the conditional distribution of $X_{t+k}$ given $X_{t}=n$ is the
convolution of a Binomial($n,\alpha^{k}$) distribution and the $m$-fold
convolution of the Poissonian Binomial ($k,\alpha,p$) distribution of
(\ref{pmf_Po_Bin}).

By (\ref{mean_var}), the mean, the variance and the index of dispersion of
$X_{t}$ are
\[
\mu_{X}=\frac{mp}{1-\alpha},\quad\sigma_{X}^{2}={\frac{mp(1+\alpha
-p)}{1-\alpha^{2}}\quad}\text{and }I_{X}=1-{\frac{p}{1+\alpha}.}%
\]
Clearly $\{X_{t}\}$ is underdispersed.

The factorial moments $(\mu_{\lbrack r]},n\geq1)$ of $X_{t}$ can be obtained
from the version (\ref{Bin1}) of $\varphi_{X}(z)$. In this case,
\[
\varphi_{X}(1+t)=\Biggl[1+\sum_{r=1}^{\infty}{\frac{p^{r}\alpha^{\binom{r}{2}%
}}{\prod_{i=1}^{r}(1-\alpha^{i})}}\cdot t^{r}\Biggr]^{m},
\]
is a power of a power series. Therefore, it admits a power series expansion
whose coefficients can be determined relatively easily for small exponents of
$t^{r}$, but not for large values of $r$. These coefficients can also be
derived via recurrence formulas (see Knopp, 1990). We won't pursue this
approach. Instead, we proceed to derive simpler recurrence formulas for
$(\mu_{\lbrack r]},n\geq1)$ by using the representation (\ref{Bin2}) of
$\varphi_{X}(z)$.

Let
\begin{equation}
\phi(z)=\sum_{n=1}^{\infty}{\frac{p^{n}}{n(1-\alpha^{n})}}(1-z)^{n}.
\label{log_Bern}%
\end{equation}

The series (\ref{log_Bern}) converges uniformly over the interval $(0,1)$ due
to the fact
\[
{\frac{p^{n}}{n(1-\alpha^{n})}} (1-z)^{n} \le{\frac{p^{n}}{n(1-\alpha^{n})}}%
\]

for every $z\in(0,1)$ and that $\sum\limits_{n=0}^{\infty}{\frac{p^{n}%
}{n(1-\alpha^{n})}}$converges. It follows that $\phi^{\prime}(z)$ and
subsequent higher order derivatives exist and converge uniformly over $(0,1)$
(see Knopp. 1990). The $r$-th derivative of $\phi(z)$ admits the representation%

\begin{equation}
\phi^{(r)}(z)=(-1)^{r}\sum_{n=r}^{\infty}{\frac{p^{n}}{1-\alpha^{n}}}%
{\frac{(n-1)!}{(n-r)!}}(1-z)^{n-r}\qquad(r\geq1). \label{derive_phi}%
\end{equation}

Uniform convergence allows for the interchange of limit (as $z\uparrow1$) and
summation in (\ref{derive_phi}). Hence,%

\begin{equation}
\phi^{(r)}(1)=(-1)^{r}{\frac{(r-1)!p^{r}}{1-\alpha^{r}}}\qquad(r\geq1).
\label{phi at 1}%
\end{equation}

Since $\ln\varphi_{X}(z)=-m\phi(z)$ (by (\ref{Bin2})), it follows that
$\varphi_{X}^{\prime}(z)=-m\varphi_{X}(z)\phi^{\prime}(z)$. An induction
argument shows that the $r^{th}$ derivative, $\varphi_{X}^{(r)}(z),$ of
$\varphi_{X}(z)$ can be obtained by the following forward recursion (with
$\varphi_{X}^{(0)}(z)=\varphi_{X}(z)$ and $\binom{0}{0}=1)$:%

\begin{equation}
\varphi_{X}^{(r)}(z)=-m\sum_{j=0}^{r-1}\binom{r-1}{j}\varphi_{X}^{(j)}%
(z)\phi^{(r-j)}(z). \label{derive_varphi}%
\end{equation}

Therefore, the factorial moments $\mu_{\lbrack r]}=\varphi_{X}^{(r)}(1)$,
$r\geq1$, are finite and satisfy the recurrence relation (with $\mu
_{\lbrack0]}=1$),
\begin{equation}
\mu_{\lbrack r]}=-m\sum_{j=0}^{r-1}\binom{r-1}{j}\mu_{\lbrack j]}\phi
^{(r-j)}(1) \qquad(r\geq1). \label{fact_mom_Bin}%
\end{equation}

For example, we have
\begin{equation}
\mu_{\lbrack1]}={\frac{mp}{1-\alpha}}\quad\text{and}\quad\mu_{\lbrack
2]}={\frac{m[(m+1)\alpha+m-1]p^{2}}{(1-\alpha)(1-\alpha^{2})}}.
\label{facmom_Bin}%
\end{equation}

Similarly to the Bernoulli case, the moments of $X_{t}$, $\mu_{r}=E(X_{t}%
^{r})$, $r\geq1$, are finite and can be obtained from their factorial
counterparts via (\ref{momfacmom}) and (\ref{Stirling_2}).

By (\ref{Bin2}), the fcgf of $X_{t}$ is
\[
\ln\varphi_{X}(1+t)=m\sum_{r=1}^{\infty}{\frac{(-1)^{r+1}(r-1)!p^{r}%
}.{(1-\alpha^{r})}}\cdot{\frac{t^{r}}{r!}}.
\]

which leads to the following formula for the factorial cumulants
$(\kappa_{\lbrack r]},r\geq1)$ of $\{X_{t}\}$:%

\begin{equation}
\kappa_{\lbrack r]}=(-1)^{r+1}{\frac{m(r-1)!p^{r}}{(1-\alpha^{r})}}
\qquad(r\geq1). \label{fact_cum_Bin}%
\end{equation}

The cumulants of $X_{t}$, $(\kappa_{r},r\geq1)$ are given by (see
(\ref{cumfaccum_1}) and (\ref{Stirling_2}))
\begin{equation}
\kappa_{r}=m\sum_{j=0}^{r}S(r,j)(-1)^{j+1}{\frac{(j-1)!p^{j}}{(1-\alpha^{j})}
\qquad(r\ge1) \label{cum_Bin}}%
\end{equation}

\section{Stationary \inar1 models with\textbf{ } Poissonian Binomial\textbf{
}innovations}

In this section, we develop a stationary \inar1 process with a Poissonian
Binomial innovation sequence with pgf and pmf given in (\ref{pgf_Po_Bin}) and
\ref{pmf_Po_Bin}), respectively. It is a generalization of the stationary
\inar1 process with binomial innovations seen in the previous Section.

\textbf{Theorem 7. }Let $\{X_{t}\}$ be the stationary \inar1 process driven by
(\ref{inar1_eq}) and with a Poissonian Binomial($m,q,c$) innovation sequence
for some positive integer $m$ and some real numbers $q,c \in(0,1)$. Then the
marginal pgf $\varphi_{X}(z)$ of $\{X_{t}\}$ admits the following representations:

\begin{enumerate}
\item
\begin{equation}
\varphi_{X}(z)=\prod_{j=0}^{m-1}\Bigl[1+\sum_{n=1}^{\infty}{\frac{(cq^{j})^{n}
(z-1)^{n}\alpha^{\binom{n}{2}}}{\prod_{l=1}^{n}(1-\alpha^{l})}}\Bigr].
\label{Po_Bin1}%
\end{equation}

\item
\begin{equation}
\varphi_{X}(z)=\exp\left\{  -\sum_{n=1}^{\infty}{\frac{1-q^{mn}}{1-q^{n}}%
}{\frac{c^{n}}{{n(1-\alpha^{n}})}}(1-z)^{n}\right\}  . \label{Po_Bin2}%
\end{equation}

\end{enumerate}

\textbf{Proof:} Let $\Psi(z)$ be the pgf of the Poissonian Binomial ($m,q,c)$
distribution as given in (\ref{pgf_Po_Bin}). Then,
\[
\Psi(1-\alpha^{i}+\alpha^{i}z)=\prod_{j=0}^{m-1}(1+c\alpha^{i}q^{j}(z-1)),
\]
which is the pgf of a Poissonian Binomial($m,q,c\alpha^{i}$). By Theorem 1,
the marginal pgf $\varphi_{X}(z)$ is
\[
\varphi_{X}(z)=\prod_{j=0}^{m-1}\prod_{i=0}^{\infty}(1+c\alpha^{i}%
q^{j}(z-1)).
\]
Noting that $\prod\limits_{i=0}^{\infty}(1+c\alpha^{i}q^{j}(z-1))$ is the
marginal pgf of a stationary \inar1 process with Bernoulli($cq^{j}$)
innovations, representations (\ref{Po_Bin1}) and (\ref{Po_Bin2}) follow from
(\ref{Bern1}) and (\ref{Bern2}), respectively. $\hfill\blacksquare$

For each $j\ge0$, we denote by $\{q_{r}^{(j)}\}$ the pmf with pgf
\[
\varphi_{j}(z)=1+\sum_{n=1}^{\infty}{\frac{(cq^{j})^{n} (z-1)^{n}%
\alpha^{\binom{n}{2}}}{\prod_{l=1}^{n}(1-\alpha^{l})}}.
\]

By Theorem 5 and (\ref{pmf_Bern}),
\begin{equation}
q_{r}^{(j)}=\sum_{k=r}^{\infty}(-1)^{k-r}\binom{k}{r}{\frac{(cq^{j})^{k}%
\alpha^{\binom{k}{2}}}{\prod_{l=1}^{k}(1-\alpha^{l})},}\quad r\geq0.
\label{pmf_Bernj}%
\end{equation}

It follows by Theorem 7 and (\ref{Po_Bin1}) that the marginal pmf $\{q_{r}\}$
of the stationary \inar1 process with a Poisson Binomial($m,q,c$) innovation
results from the convolution of the pmf's $(\{q_{r}^{(j)}\},0\leq j\leq m-1)$,
i.e.,
\begin{equation}
q_{r}=(q^{(0)}\ast q^{(1)}\ast\cdots\ast q^{(m-1)})_{r}\quad(r\geq0).
\label{pmf_Bo_Bin_inar}%
\end{equation}

We proceed to give additional properties for the marginal distribution of the
stationary \inar1 process $\{X_{t}\}$ with Poisson Binomial($m,q,c$) innovations.

By (\ref{MarCh}), the 1-step transition probability is given by%

\begin{equation}
P(X_{t}=k|X_{t-1}=l)=\sum_{j=\max(k-m,0)}^{\min(l,k)}\binom{l}{j}\alpha^{j}
(1-\alpha)^{l-j} q_{k-j}(m,q,c). \label{Tran_Po_Bin}%
\end{equation}

By (\ref{k_step_pgf}), the conditional pgf of $X_{t+k}$ given $X_{t}$
satisfies%
\[
\varphi_{X_{t+k}|X_{t}}(z)=\left(  1-\alpha^{k}+\alpha^{k}z\right)  ^{X_{t}%
}\times\prod_{j=0}^{m-1}\Bigl[\prod_{i=0}^{k-1}(1-(cq^{j})\alpha
^{i}(1-z))\Bigr].
\]
Therefore, the conditional distribution of $X_{t+k}$ given $X_{t}=n$ is the
convolution of a Binomial($n,\alpha^{k}$) distribution and the Poissonian
Binomial ($k, \alpha, cq^{j}$) distributions, $j=0,1,\cdots,m-1$.

The innovation sequence of $\{X_{t}\}$ being Poissonian Binomial($m,q,c)$ is
underdispersed with mean and variance (see Kemp, 1987) given by
\[
\mu_{\varepsilon}={\frac{(1-q^{n})c}{1-q}}\quad\hbox{and}\quad\sigma
_{\varepsilon}^{2}={\frac{(1-q^{n})c}{1-q}}-{\frac{(1-q^{2n})c^{2}}{1-q^{2}}%
}.
\]
Therefore, the marginal distribution of $\{X_{t}\}$ is also underdispersed
with mean and variance and dispersion index given by (\ref{mean_var}):
\[
\mu_{X}={\frac{(1-q^{n})c}{(1-\alpha)(1-q)}},\quad\sigma_{X}^{2}%
={\frac{(1-q^{n})c}{(1-\alpha)(1-q)}}-{\frac{(1-q^{2n})c^{2}}{(1-\alpha
^{2})(1-q^{2})}}%
\]
and
\[
I_{X}=1-{\frac{(1+q^{n})c}{(1+\alpha)(1+q)}.}%
\]

The factorial moments $(\mu_{\lbrack r]},n\geq1)$ of $X_{t}$ can be obtained
from the version (\ref{Po_Bin1}) of $\varphi_{X}(z)$. In this case,
\[
\varphi_{X}(1+t)=\prod_{j=0}^{m-1}\Bigl[1+\sum_{n=1}^{\infty}{\frac
{(cq^{j})^{n} \alpha^{\binom{n}{2}}}{\prod_{l=1}^{n}(1-\alpha^{l})}t^{n}%
}\Bigr].
\]
is a finite product of a power series. Therefore, it admits a power series
expansion whose coefficients have rather complex expressions, even via
recurrence formulas.

We proceed to derive simpler recurrence formulas for $(\mu_{\lbrack r]}%
,n\geq1)$ by using the representation (\ref{Po_Bin2}) of $\varphi_{X}(z)$.
Let
\begin{equation}
\phi_{1}(z) =-\ln\varphi_{X}(z)=\sum_{n=1}^{\infty}{\frac{1-q^{mn}}{1-q^{n}}%
}{\frac{c^{n}}{n(1-\alpha^{n})}}(1-z)^{n}. \label{log_Po_Bin}%
\end{equation}

The argument we used to derive (\ref{derive_phi})--(\ref{fact_mom_Bin}) in
Section 5 carries over almost verbatim. We state the main results without
further explanations.

The $r$-th derivative of $\phi_{1}(z)$ admits the representation%

\begin{equation}
\phi_{1}^{(r)}(z)=(-1)^{r}\sum_{n=r}^{\infty}{\frac{(1-q^{mn})c^{n}}%
{(1-q^{n})(1-\alpha^{n})}}{\frac{(n-1)!}{(n-r)!}}(1-z)^{n-r}\qquad r\geq1.
\label{derive_phi_1}%
\end{equation}

Hence,%

\begin{equation}
\phi_{1}^{(r)}(1)=(-1)^{r}{\frac{(1-q^{mr})c^{r} }{(1-q^{r}) (1-\alpha^{r})}%
}(r-1)!\qquad(r\geq1). \label{phi_1 at 1}%
\end{equation}

Since $\ln\varphi_{X}(z)=-\phi_{(}z)$, the $r^{th}$ derivative, $\varphi
_{X}^{(r)}(z),$ of $\varphi_{X}(z)$ can be obtained by the following forward
recursion (with $\varphi_{X}^{(0)}(z)=\varphi_{X}(z)$ and $\binom{0}{0}=1)$:
\begin{equation}
\varphi_{X}^{(r)}(z)=-\sum_{j=0}^{r-1}\binom{r-1}{j}\varphi_{X}^{(j)}%
(z)\phi_{1}^{(r-j)}(z). \label{derive-varphi1}%
\end{equation}

The factorial moments $\mu_{\lbrack r]}=\varphi_{X}^{(r)}(1)$, $r\geq1$, are
finite and satisfy the recurrence relation (with $\mu_{\lbrack0]}=1$),
\begin{equation}
\mu_{\lbrack r]}=-\sum_{j=0}^{r-1}\binom{r-1}{j}\mu_{\lbrack j]}\phi_{1}
^{(r-j)}(1) \qquad(r\geq1). \label{fact_mom_Po_Bin}%
\end{equation}

The moments of $X_{t}$, $\mu_{r}=E(X_{t}^{r})$, $r\geq1$, are finite and can
be obtained from their factorial counterparts via (\ref{momfacmom}) and
(\ref{Stirling_2}).

We see by (\ref{Po_Bin2}) that the fcgf of $X_{t}$ is
\[
\ln\varphi_{X}(1+t)=\sum_{n=1}^{\infty}{\frac{(-1)^{n+1}(n-1)!(1-q^{mn})c^{n}%
}{(1-q^{n})(1-\alpha^{n})}}{\frac{t^{n}}{n!}},
\]

which yields the following formula for the factorial cumulants $(\kappa
_{\lbrack r]},r\geq1)$ of $\{X_{t}\}$:
\begin{equation}
\kappa_{\lbrack r]}={\frac{(-1)^{r+1}(r-1)!(1-q^{mr})c^{r} }{(1-q^{r}%
)(1-\alpha^{r})}}\qquad(r\geq1). \label{fact_cum_Po_Bin}%
\end{equation}

The cumulants of $X_{t}$, $(\kappa_{r},r\geq1)$ are given by (see
(\ref{cumfaccum_1}) and (\ref{Stirling_2}))
\begin{equation}
\kappa_{r}=\sum_{j=0}^{r}S(r,j){\frac{(-1)^{j+1}(j-1)!(1-q^{mj})c^{j}%
}{(1-q^{j})(1-\alpha^{j})},}\qquad r\geq1. \label{cum_P0_Bin}%
\end{equation}

\section{Stationary \inar1 processes with\textbf{ } Heine\textbf{
}innovations}

A distribution on $\mathbb{Z}_{+}$ is said to have the Heine distribution
(Heine($\lambda,q$)) with parameters $\lambda>0$ and $q\in(0,1)$ if its pgf
and pmf are respectively given by
\begin{equation}
\Psi(z)=\prod_{j=0}^{\infty}\bigl(1-\beta_{j}+\beta_{j}%
z\bigr) \label{Heine_pgf}%
\end{equation}
and%
\begin{equation}
f_{r}={\frac{\lambda^{r}q^{r(r-1)/2}}{\prod_{l=1}^{r}(1-q^{l})}}%
f_{0},r\geq1\quad\text{and }f_{0}=\prod_{j=0}^{\infty}(1-\lambda q^{j})^{-1},
\label{Heine_pmf}%
\end{equation}
where $\beta_{j}={\frac{\lambda q^{j}}{1+\lambda q^{j}}}$ for $j\geq0$.

The Heine distribution was introduced by Benkherouf and Bather (1988). Kemp
(1992) studied many of its properties. More details can be found in these
references and in Johnson et al. (2005, Section 10.8.2). The Heine
distribution is underdispersed and its mean and variance are
\begin{equation}
\mu=\sum_{r=0}^{\infty}{\frac{\lambda q^{r}}{1+\lambda q^{r}}}\quad
\hbox{and}\quad\sigma^{2}=\sum_{r=0}^{\infty}{\frac{\lambda q^{r}}{(1+\lambda
q^{r})^{2}}}. \label{Heine_mean_var}%
\end{equation}

We recall a few facts about double infinite products. Let $\{a_{mn}\}$ be a
double sequence. The double infinite product $\prod\limits_{i=0}^{\infty}%
\prod\limits_{j=0}^{\infty}(1+a_{ij})$ is defined as the limit of the double
sequence $P_{mn}=\prod\limits_{i=0}^{m}\prod\limits_{j=0}^{n}(1+a_{ij})$ as
$m,n\rightarrow\infty$. If $\sum\limits_{i=0}^{\infty}\sum\limits_{j=0}%
^{\infty}|a_{ij})|<\infty$, then the double infinite product $\prod
\limits_{i=0}^{\infty}\prod\limits_{j=0}^{\infty}(1+a_{ij})$ converges.
Moreover, if $\prod\limits_{i=0}^{\infty}\prod\limits_{j=0}^{\infty}%
(1+a_{mn})$ and the iterated infinite products
\[
\prod\limits_{i=0}^{\infty}\Bigl[\prod\limits_{j=0}^{\infty}\bigl(1+a_{ij}%
\bigr)\Bigr]\quad\hbox{and}\quad\prod\limits_{j=0}^{\infty}\Bigl[\prod
\limits_{i=0}^{\infty}\bigl(1+a_{ij}\bigr)\Bigr]
\]
converge, then they converge to the same limit.

\textbf{Theorem 8. }Let $\{X_{t}\}$ be the stationary \inar1 process driven by
(\ref{inar1_eq}) and with a Heine($\lambda,q)$ innovation sequence for some
$\lambda>0$ and $0<q<1$. Then the marginal pgf $\varphi_{X}(z)$ of $\{X_{t}\}$
admits the following representations:

\begin{enumerate}
\item
\begin{equation}
\varphi_{X}(z)=\prod_{j=0}^{\infty}\Bigl[1+\sum_{n=1}^{\infty}{\frac{\beta
_{j}^{n}(z-1)^{n}\alpha^{\binom{n}{2}}}{\prod_{l=1}^{n}(1-\alpha^{l})}}\Bigr],
\label{Heine1}%
\end{equation}
where $\beta_{j}$ is as in (\ref{Heine_pgf}).

\item
\begin{equation}
\varphi_{X}(z)=\exp\left\{  -\sum_{n=1}^{\infty}{\frac{B_{n}}{n(1-\alpha^{n}%
)}}(1-z)^{n}\right\}  \label{Heine2}%
\end{equation}
with $B_{n}=\sum_{j=0}^{\infty}\beta_{j}^{n}$, $n\geq1$.
\end{enumerate}

\textbf{Proof:} First, we note that $0<\beta_{j}<1$ for any $j\geq0$.
Moreover, for any $n\geq1$,
\begin{equation}
B_{n}=\sum_{j=0}^{\infty}{\frac{\lambda^{n}(q^{n})^{j}}{(1+\lambda q^{j})^{n}%
}}\leq{\frac{\lambda^{n}}{1-q^{n}}}<\infty. \label{B_inequal}%
\end{equation}
The pgf $\Psi(z)$ of the innovation sequence of $\{X_{t}\}$ is given by the
right hand side of (\ref{Heine_pgf}). Noting that for every $i\geq0$,
\[
\Psi(1-\alpha^{i}+\alpha^{i}z)=\prod_{j=0}^{\infty}\bigl(1-\beta_{j}\alpha
^{i}(1-z))\bigr),
\]
it follows by Theorem 3 and (\ref{varphi_thm3}) that
\begin{equation}
\varphi_{X}(z)=\prod_{i=0}^{\infty}\Bigl[\prod_{j=0}^{\infty}\bigl(1-\beta
_{j}\alpha^{i}(1-z)\bigr)\Bigr]. \label{Heine3}%
\end{equation}
Clearly, the right hand side of (\ref{Heine3}) converges. A straightforward
argument shows that the double infinite product $\prod\limits_{i=0}^{\infty
}\prod\limits_{j=0}^{\infty}\bigl(1-\beta_{j}\alpha^{i}(1-z)\bigr)$ converges.
In order to be able to interchange the order of the infinite products in
(\ref{Heine3}), it remains to show that the iterated infinite product
\[
\prod\limits_{j=0}^{\infty}\Bigl[\prod\limits_{i=0}^{\infty}\bigl(1-\beta
_{j}\alpha^{i}(1-z)\bigr)\Bigr]=\prod\limits_{j=0}^{\infty}P_{j}(z)
\]
converges, where for each $j\geq0$,%
\[
P_{j}(z)=\prod\limits_{i=0}^{\infty}\bigl(1-\beta_{j}\alpha^{i}(1-z)\bigr).
\]
Note that for each $j\geq0,$ $P_{j}(\cdot)$ has the form of the pgf of the
marginal of an \inar1 process with a Bernoulli($\beta_{j}$) innovation (see
(\ref{Bern3}) in Appendix). Therefore, by Theorem 5 and (\ref{Bern1}),
\begin{equation}
P_{j}(z)=1+\sum_{n=1}^{\infty}{\frac{\beta_{j}^{n}(z-1)^{n}\alpha^{\binom
{n}{2}}}{\prod_{l=1}^{n}(1-\alpha^{l})}}. \label{Heine}%
\end{equation}
For $j\geq0$, denote
\[
\zeta_{j}(z)=\sum_{n=1}^{\infty}{\frac{\beta_{j}^{n}(z-1)^{n}\alpha^{\binom
{n}{2}}}{\prod_{l=1}^{n}(1-\alpha^{l})}}.
\]
By (\ref{B_inequal}) and $0\leq z\leq1$, we have
\[
\sum\limits_{j=0}^{\infty}|\zeta_{j}(z)|\leq\sum_{n=1}^{\infty}\sum
_{j=0}^{\infty}{\frac{\beta_{j}^{n}\alpha^{\binom{n}{2}}}{\prod_{l=1}%
^{n}(1-\alpha^{l})}}=\sum_{n=1}^{\infty}{\frac{B_{n}\alpha^{\binom{n}{2}}%
}{\prod_{l=1}^{n}(1-\alpha^{l})}}\leq\sum_{n=1}^{\infty}a_{n},
\]
where
\[
a_{n}={\frac{\lambda^{n}}{1-q^{n}}}{\frac{\alpha^{\binom{n}{2}}}{\prod
_{l=1}^{n}(1-\alpha^{l})}.}%
\]
Recall that $\alpha,q\in(0,1)$ and hence
\[
\lim_{n\rightarrow\infty}{\frac{a_{n+1}}{a_{n}}}=\lim_{n\rightarrow\infty
}{\frac{\lambda(1-q^{n})\alpha^{n}}{(1-q^{n+1})(1-\alpha^{n+1})}}=0.
\]
By the ratio test, $\sum\limits_{j=0}^{\infty}|\zeta_{j}(z)|$ converges
uniformly over $z\in\lbrack0,1]$. This in turn implies (see Knopp, 1990) that
$\prod\limits_{j=0}^{\infty}P_{j}(z)$ converges uniformly over $z\in
\lbrack0,1]$. The representation (\ref{Heine1}) then follows by interchanging
the order of the infinite products in (\ref{Heine3}) and using (\ref{Heine}).
We now prove (\ref{Heine2}). By the first part of the proof, we have
\[
\varphi_{X}(z)=\prod_{j=0}^{\infty}\Bigl[\prod_{i=0}^{\infty}\bigl(1-\beta
_{j}\alpha^{i}(1-z)\bigr)\Bigr].
\]
Applying the representation (\ref{Bern2}) to $\prod_{i=0}^{\infty
}\bigl(1-\beta_{j}\alpha^{i}(1-z)\bigr)$ with $p=\beta_{j}$, we have
\begin{equation}
\varphi_{X}(z)=\exp\left\{  -\sum_{j=0}^{\infty}\Bigl[\sum_{n=1}^{\infty
}{\frac{\beta_{j}^{n}}{n(1-\alpha^{n})}}(1-z)^{n}\Bigr]\right\}  .
\label{Heine5}%
\end{equation}
This implies that the double series in (\ref{Heine5}) is convergent. Since its
terms are nonnegative (as $0\leq z\leq1$), the order of summation can be
interchanged (by Cauchy's criterion for double series). This establishes the
representation (\ref{Heine2}).$\hfill\blacksquare$

For each $j\ge0$, we denote by $\{p_{r}^{(j)}\}$ the pmf with pgf
\[
\varphi_{j}(z)=1+\sum\limits_{n=1}^{\infty}{\frac{\beta_{j}^{n}(z-1)^{n}
\alpha^{\binom{n}{2}}}{\prod_{l=1}^{n} (1-\alpha^{l})}}.
\]

By Theorem 5 and (\ref{pmf_Bern}),
\begin{equation}
p_{r}^{(j)}=\sum_{k=r}^{\infty}(-1)^{k-r}\binom{k}{r}{\frac{\beta_{j}%
^{k}\alpha^{\binom{k}{2}}}{\prod_{l=1}^{k}(1-\alpha^{l})},}\quad
r\geq0.\label{pmf_Bernj1}%
\end{equation}

It follows by Theorem 8, (\ref{Heine1}) and the Continuity Theorem that the
marginal pmf $\{p_{r}\}$ of the stationary \inar1 process with a
Heine($\lambda,q$) innovation is
\begin{equation}
p_{r}=\lim_{j\rightarrow\infty}(p^{(0)}\ast p^{(1)}\ast\cdots\ast p^{(j)}%
)_{r}\quad(r\ge0). \label{pmf_Heine_inar}%
\end{equation}

Next, we discuss several properties of the marginal distribution of the
stationary \inar1 process $\{X_{t}\}$ with Heine$(\lambda,q)$ innovations.

The 1-step transition probability can be obtained from (\ref{MarCh}). Given
there are no notable simplifications of the formulas, we omit the details.

By (\ref{k_step_pgf}), the conditional distribution of $X_{t+k}$ given
$X_{t}=n$ results from the convolution of of $k+1$ distributions, namely a
Binomial($n,\alpha^{k}$) distribution and the distributions $(\{g_{r}%
^{(i)}\},0\leq i\leq k-1)$ defined as follows:
\[
g_{r}^{(i)}=\alpha^{ir}\sum_{l=0}^{\infty}\binom{r+l}{r}(1-\alpha)^{l}%
f_{r+l},
\]
where $\{f_{r}\}$ is the pmf of the Heine($\lambda,q$) distribution
(\ref{Heine_pmf}).

By (\ref{mean_var}) and (\ref{Heine_mean_var}), the mean, the variance and the
index of dispersion of $X_{t}$ are given by
\[
\mu_{X}={\frac{1}{1-\alpha}}\sum_{r=0}^{\infty}{\frac{\lambda q^{r}}{1+\lambda
q^{r}},}\quad\sigma_{X}^{2}={\frac{1}{1-\alpha^{2}}}\sum_{r=0}^{\infty}%
{\frac{\lambda q^{r}}{(1+\lambda q^{r})^{2}}}%
\]
and
\[
I_{X}=\sum\limits_{r=0}^{\infty}{\frac{\lambda q^{r}}{1+\lambda q^{r}}%
}([1+\lambda q^{r}]^{-1}+\alpha)\Big /\Bigl[(1+\alpha)\sum\limits_{r=0}%
^{\infty}{\frac{\lambda q^{r}}{1+\lambda q^{r}}}\Bigr]
\]

Since the Heine distribution is underdispersed, the \inar1 process with a
Heine innovation is underdispersed.

The factorial cumulants $(\kappa_{\lbrack r]},r\geq1)$ of $X_{t}$ are easily
obtained. Indeed, by (\ref{Heine2}), the fcgf of $X_{t}$ is
\begin{equation}
\ln\varphi_{X}(1+t)=\sum_{r=1}^{\infty}{\frac{(-1)^{r+1}(r-1)!B_{r}}%
{1-\alpha^{r}}}{\frac{t^{r}}{r!}.} \label{fcgf_Heine}%
\end{equation}
Since the series above converges everywhere, the factorial cumulants of all
orders are finite and are given by
\begin{equation}
\kappa_{\lbrack r]}=(-1)^{r+1}{\frac{(r-1)!B_{r}}{1-\alpha^{r}}}\qquad
(r\geq1). \label{fact_cum_Heine}%
\end{equation}

The cumulants of $X_{t}$, $(\kappa_{r},r\geq1)$, can be obtained from the
factorial cumulants via (\ref{cumfaccum_1}) and (\ref{Stirling_2})
\begin{equation}
\kappa_{r}=\sum_{j=0}^{r}S(r,j)(-1)^{j+1}{\frac{(j-1)!p^{j}}{(1-\alpha^{j})}}
\qquad(r\geq1) \label{cumfaccum}%
\end{equation}

The moments $\mu_{r},r\geq1$ can be computed using the formula due to Smith
(1995):
\begin{equation}
\mu_{r}=\sum_{j=0}^{r-1}\binom{r-1}{j}\kappa_{r-j}\mu_{j} \label{mom_cum}%
\end{equation}
with initial conditions $\mu_{0}=1$ and $\mu_{1}=\kappa_{1}$. In turn, the
factorial moments $\mu_{\lbrack r]},r\geq1$, of $X_{t}$ can be obtained via
the formula (see Johnson et al. (2005), Section 1.2.7):
\begin{equation}
\mu_{\lbrack r]}=\sum_{j=0}^{r}s(r,j)\mu_{j}, \label{facmom_mom}%
\end{equation}
where $\{s(r,j)\}$ are the Stirling numbers of the first kind of
(\ref{Stirling_1}).

\section{Extensions via convolution}

Let $\nu\geq1$ be a positive integer. Assume that $\Psi(z)$ in Theorem 1 is
the pgf of a finite convolution of $\nu$ pmf's $\{f_{1r}\},\{f_{2r}%
\},\cdots,\{f_{\nu r}\}$ with respective pgf's $\Psi_{1}(z),\Psi_{2}%
(z),\cdots,\Psi_{\nu}(z)$ and $\Psi_{k}^{\prime}(1)<\infty$ for every $1\leq
k\leq\nu$. We denote by $\{f_{kr}^{(i)}\}$ the pmf with pgf $\Psi_{k}%
(1-\alpha^{i}+\alpha^{i}z)$ (with $\{f_{kr}^{(0)}\}=\{f_{kr}\}$). Note that
the pmf $\{f_{r}^{(i)}\}$ of $\Psi(1-\alpha^{i}+\alpha^{i}z)$ is the
convolution
\begin{equation}
f_{r}^{(i)}=\Bigl(f_{1\cdot}^{(i)}\ast f_{2\cdot}^{(i)}\ast\cdots\ast
f_{\nu\cdot}^{(i)}\Bigr)_{r},\quad r\geq0. \label{pmf_conv_fri}%
\end{equation}
Applying Theorem 1 to $\Psi_{k}(z)$ for each $k\in\{1,2,\cdots,\nu\}$, the
function
\begin{equation}
\varphi^{(k)}(z)=\prod\limits_{i=0}^{\infty}\Psi_{k}(1-\alpha^{i}+\alpha^{i}z)
\label{varphi_k}%
\end{equation}
is a pgf and its pmf is
\begin{equation}
p_{r}^{(k)}=\lim_{n\rightarrow\infty}\Bigl(f_{k\cdot}^{(0)}\ast f_{k\cdot
}^{(1)}\ast\cdots\ast f_{k\cdot}^{(n-1)}\Bigr)_{r}\quad(r\geq0),
\label{pmf_varphi_k}%
\end{equation}
The pgf $\varphi(z)$ of (\ref{varphi_thm1}) with $\Psi(z)=\prod_{k=1}^{\nu
}\Psi_{k}(z)$ is $\varphi(z)=\prod\limits_{k=1}^{\nu}\varphi^{(k)}(z)$.
Therefore, its pmf
\begin{equation}
p_{r}=\Bigl(p^{(1)}\ast p^{(2)}\ast\cdots\ast p^{(\nu)}\Bigr)_{r}
\label{pmf_conv}%
\end{equation}
results from the convolution of $(\{p_{r}^{(k)}\},1\leq k\leq\nu).$

The mean, the variance, the cumulants and the factorial cumulants of
$\{p_{r}\}$ are the sums of their $(\{p_{r}^{(k)}\},1\leq k\leq\nu)$
counterparts assuming the latter are finite. The moments and the factorial
moments could possibly be computed using the equations (\ref{mom_cum}) and
(\ref{facmom_mom}).

The Binomial and Poissonian distributions being finite convolutions of Bernoulli 
distributions, the \inar1 processes introduced in Sections 5 and 6 could have been 
developed using the approach described above in conjonction with the Bernoulli
\inar1 of Section 3. However, the authors deemed  these two models important 
enough to be treated separately and more thoroughly.  

Next, we discuss some simple examples of stationary \inar1 processes whose
innovation is the convolution of a Poisson($\lambda$) and the underdispersed
distributions discussed in Sections 3-7. In the enumeration that follows,
$\{X_{t}\}$ designates a stationary \inar1 process of (\ref{inar1_eq}).

\begin{enumerate}
\item If the innovation $\{\varepsilon_{t}\}$ admits marginal law
$\varepsilon_{t}\sim Pois(\lambda)\ast Logarithmic(p)$, $0<p<1-1/e$, then its
marginal distribution will result from the convolution of a Poisson(${\frac
{\lambda}{1-\alpha}}$) and the pmf $\{p_{r}\}$ in Theorem 4.

\item If the innovation $\{\varepsilon_{t}\}$ has the Power-Law distribution
of the first kind ($PL_{1}(\lambda,p)$), i.e., $\varepsilon_{t}\sim
Pois(\lambda)\ast Bernoulli(p)$, $0<p<1$, then its marginal distribution will
result from the convolution of a Poisson(${\frac{\lambda}{1-\alpha}}$) and the
pmf $\{p_{r}\}$ of (\ref{pmf_Bern}) in Theorem 5. Note that this model was
applied in Section 2.3 of Wei{\ss } (2013).

\item If the innovation $\{\varepsilon_{t}\}$ has a Power-Law distribution of
order $m$ ($PL_{m}^{\ast}((\lambda,p)$), i.e., $\varepsilon_{t}\sim
Pois(\lambda)\ast Binomial(m,p)$, then the marginal of the corresponding
stationary \inar1 process is the convolution of Poisson(${\frac{\lambda
}{1-\alpha}}$) and the pmf $\{p_{r}\}$ of (\ref{pmf_Bin}) in Theorem 6.

\item If the innovation $\{\varepsilon_{t}\}$ admits marginal law
$\varepsilon_{t}\sim Pois(\lambda_{1})\ast Heine(\lambda,q)$, $\lambda>0$ and
$0<q<1$, then its marginal distribution will result from the convolution of a
Poisson(${\frac{\lambda_{1}}{1-\alpha}}$) and the pmf $\{p_{r}\}$ of
(\ref{pmf_Heine_inar}) of Theorem 8.
\end{enumerate}

\vskip .5 cm

\centerline{\bf Appendix}

This section is devoted to the proof of Theorem 5. We start out with a Lemma.

\textbf{Lemma 2}

Assume $n\geq2$ and $a_{i}\in(0,1)$ for $i=0,1,2,\cdots,n-1$. Then,

\begin{enumerate}
\item
\begin{equation}
\prod_{i=0}^{n-1}(1-a_{i})=1+\sum_{k=1}^{n}(-1)^{k}\Bigl[\sum\limits_{0\leq
j_{1}<j_{2}<\cdots<j_{k}\leq n-1}\prod_{l=1}^{k}a_{j_{l}}\Big]. \label{temp_1}%
\end{equation}

\item
\begin{equation}
\sum_{0\leq j_{1}<j_{2}<\cdots<j_{k}\leq n-1}\alpha^{j_{1}}\alpha^{j_{2}%
}\cdots\alpha^{j_{k}}=\alpha^{\binom{k}{2}}\prod_{l=0}^{k-1}{\frac
{1-\alpha^{n-l}}{1-\alpha^{l+1}}} \label{temp_2}%
\end{equation}
for every $k\in\{1.\cdots,n\}$.
\end{enumerate}

\textbf{Proof:} (1) follows by a straightforward induction.

(2) We also proceed by induction. The result is trivially true for $n=2$
(forces $k=1$). Assume the assertion is true up to $n$. It is clear that
(\ref{temp_2}) holds for $n+1$ and $k=n+1$. As in this case,
\[
\sum_{0\leq j_{1}<j_{2}<\cdots<j_{n+1}\leq n}\alpha^{j_{1}}\alpha^{j_{2}%
}\cdots\alpha^{j_{k}}=\alpha^{\sum_{k=0}^{n}k}=\alpha^{\binom{n+1}{2}}%
=\alpha^{\binom{n+1}{2}}\prod_{l=0}^{n}{\frac{1-\alpha^{n+1-l}}{1-\alpha
^{l+1}}}.
\]
Assume now $k\in\{1,2,\cdots,n\}$. Setting $J=(j_{1},j_{2},\cdots,j_{k}%
)\in\mathbb{N}^{k}$, it is clear that
\[
\{J\in\mathbb{N}^{k}:0\leq j_{1}<j_{2}<\cdots<j_{k}\leq n\}=A\cup B,
\]
where $A=\{J\in\mathbb{N}^{k}:0\leq j_{1}<j_{2}<\cdots<j_{k}\leq n-1\}$ and
$B=\{J\in\mathbb{N}^{k}:0\leq j_{1}<j_{2}<\cdots<j_{k-1}\leq n-1,j_{k}=n\}$.
Therefore,
\[
\sum\limits_{0\leq j_{1}<j_{2}<\cdots<j_{k}\leq n}\prod_{l=1}^{k}\alpha
^{j_{l}}=\sum_{J\in A}\prod_{l=1}^{k}\alpha^{j_{l}}+\sum_{J\in B}\alpha
^{n}\prod_{l=1}^{k-1}\alpha^{j_{l}}.
\]
Using the induction hypothesis, it follows that
\[
\sum\limits_{J\in A}\prod\limits_{l=1}^{k}\alpha^{j_{l}}=\alpha^{\binom{k}{2}%
}\prod\limits_{l=0}^{k-1}{\frac{1-\alpha^{n-l}}{1-\alpha^{l+1}}}%
\]
and
\[
\sum\limits_{J\in B}\alpha^{n}\prod\limits_{l=1}^{k-1}\alpha^{j_{l}}%
=\alpha^{n}\alpha^{\binom{k-1}{2}}\prod\limits_{l=0}^{k-2}{\frac
{1-\alpha^{n-l}}{1-\alpha^{l+1}}}%
\]
which implies
\[
\sum\limits_{0\leq j_{1}<j_{2}<\cdots<j_{k}\leq n}\prod_{l=1}^{k}\alpha
^{j_{l}}={\frac{\prod\limits_{l=0}^{k-2}(1-\alpha^{n-l})\Bigl[(1-\alpha
^{n-k+1})\alpha^{\binom{k}{2}}+(1-\alpha^{k})\alpha^{n}\alpha^{\binom{k-1}{2}%
}\Bigr]}{\prod\limits_{l=0}^{k-1}(1-\alpha^{l+1})}}.
\]
Now, noting that $\binom{k}{2}=\binom{k-1}{2}+k-1$, it is easily seen that
\[
(1-\alpha^{n-k+1})\alpha^{\binom{k}{2}}+(1-\alpha^{k})\alpha^{n}\alpha
^{\binom{k-1}{2}}=\alpha^{\binom{k}{2}}(1-\alpha^{n+1}).
\]
Therefore, (\ref{temp_2}) holds for $n+1$. $\hfill\blacksquare$

\textbf{Proof of Theorem 5:}

Let $\{X_{t}\}$ be the stationary \inar1 process with a Bernoulli($p$)
innovation sequence. By Theorem 3 and (\ref{varphi_thm3}), its marginal pgf
is
\begin{equation}
\varphi_{X}(z)=\prod_{i=0}^{\infty}(1-p\alpha^{i}(1-z)). \label{Bern3}%
\end{equation}
Since $\varphi_{X}(z)=\displaystyle\lim_{n\rightarrow\infty}\prod_{i=0}%
^{n-1}(1-p\alpha^{i}(1-z))$, we conclude by the continuity theorem that the
marginal pmf $\{p_{r}\}$ of $\{X_{t}\}$ is the weak limit of a sequence of
Poissonian Binomial distributions of (\ref{pgf_Po_Bin}) and (\ref{pmf_Po_Bin}%
), with $m=n$, $q=\alpha$ and $c=p$. Let $r\geq0$. We define a purely atomic
measure, we denote $meas_{r}$, on $\mathbb{N}_{r}=\{r,r+1,r+2,\cdots\}$ and
its power set $\mathcal{P}(\mathbb{N}_{r})$ as follows:
\begin{equation}
meas_{r}(\{k\})={\frac{p^{k}\alpha^{\binom{k}{2}}}{\prod_{l=1}^{k}%
(1-\alpha^{l})},}\quad(k\geq r), \label{atomic_meas}%
\end{equation}
with $meas_{0}(\{0\})=1$. It is clear that $\sum_{k=r}^{\infty}meas_{r}%
(\{k\})<\infty$. Therefore, $meas_{r}$ is a finite measure. Define now the
sequence of functions $\{f_{n}(\cdot)\}$ on $\mathbb{N}_{r}$ by
\[
f_{n}(k)=%
\begin{cases}
(-1)^{k-r}\binom{k}{r}\prod\limits_{l=0}^{k-1}(1-\alpha^{n-l}) &
\mbox{if }k=r,r+1,\cdots,n\\
0 & \mbox{if }k>n.
\end{cases}
\]
Define $h(k)=\binom{k}{r}$ on $\mathbb{N}_{r}$. It is clear that
$|f_{n}(k)|\leq h(k)$ (recall $\alpha\in(0,1)$) and that $\sum\limits_{k=r}%
^{\infty}h(k)meas_{r}(\{k\})<\infty$ (by the ratio test). Moreover, for every
$k\in\mathbb{N}_{r}$,
\[
f(k)=\lim_{n\rightarrow\infty}f_{n}(k)=(-1)^{k-r}\binom{k}{r}.
\]
Rewriting $p_{r}^{(n)}$ in terms of the discrete integral of $f_{n}(k)$ on the
measure space $(\mathbb{N}_{r},\mathcal{P}(\mathbb{N}_{r}),meas_{r})$ and
calling on the dominated convergence theorem, we have
\[
p_{r}=\lim_{n\rightarrow\infty}\int_{\mathbb{N}_{r}}f_{n}(k)\,meas(dk)=\int%
_{\mathbb{N}_{r}}f(k)\,meas(dk),
\]
which is precisely (\ref{pmf_Bern}) and thus part (1) of the Theorem is
established. To show (2), note that
\[
P(X_{t}\geq r)=\sum_{j=r}^{\infty}\sum_{k=j}^{\infty}(-1)^{k-j}\binom{k}%
{j}{\frac{p^{k}\alpha^{\binom{k}{2}}}{\prod_{l=1}^{k}(1-\alpha^{l})}.}%
\]
Since the double series above converges absolutely, interchanging summations
is allowed, leading to
\[
P(X_{t}\geq r)=\sum_{k=r}^{\infty}\Biggl(\sum_{j=r}^{k}(-1)^{k-j}\binom{k}%
{j}\Biggr){\frac{p^{k}\alpha^{\binom{k}{2}}}{\prod_{l=1}^{k}(1-\alpha^{l})}.}%
\]
We have by induction on $k$ that $\sum_{j=r}^{k}(-1)^{k-j}\binom{k}%
{j}=(-1)^{k-r}\binom{k-1}{r-1}$, establishing (\ref{tail_Bern}).

For part (3), we note first that $\varphi_{X}(z)$ of (\ref{Bern3}) can be
rewritten as
\begin{equation}
\varphi_{X}(z)=\exp\left\{  \sum_{i=0}^{\infty}\ln(1-p\alpha^{i}%
(1-z))\right\}  . \label{Bern4}%
\end{equation}
The representation (\ref{Bern2}) of $\varphi_{X}(z)$ follows by way of the
power series expansion of
\[
-\ln(1-x)=\sum_{n=1}^{\infty}x^{n}/n,\quad0\leq x<1
\]
applied to $x=p\alpha^{i}(1-z)$ in (\ref{Bern4}).

To prove (\ref{Bern1}), we first note that by letting $a_{i}=p\alpha^{i}(1-z)$
in (\ref{temp_1}) and using (\ref{temp_2}), we obtain the following expression
for $\varphi_{n}(z)$ of (\ref{pgf_Po_Bin}):
\begin{equation}
\varphi_{n-1}(z)=1+\sum_{k=1}^{n}p^{k}(z-1)^{k}\alpha^{\binom{k}{2}}%
\prod_{l=0}^{k-1}{\frac{1-\alpha^{n-l}}{1-\alpha^{l+1}}} \label{temp_3}%
\end{equation}
and therefore,
\begin{equation}
\varphi_{X}(z)=\lim_{n\rightarrow\infty}\Bigl[1+\sum_{k=1}^{n}(z-1)^{k}%
\prod_{l=0}^{k-1}(1-\alpha^{n-l}){\frac{p^{k}\alpha^{\binom{k}{2}}}%
{\prod_{l=1}^{k}(1-\alpha^{l})}}\Big]. \label{temp_4}%
\end{equation}
We proceed as in the proof of (\ref{pmf_Bern}). We define a sequence of
functions $g_{n}(k)$ on the finite measure space $(\mathbb{N},\mathcal{P}%
(\mathbb{N}),meas_{0})$, where $meas_{0}$ is defined in (\ref{atomic_meas}):
\[
g_{n}(k)=%
\begin{cases}
1 & \mbox{if }k=0\\
(z-1)^{k}\prod\limits_{l=0}^{k-1}(1-\alpha^{n-l}) & \mbox{if }1\leq k\leq n\\
0 & \mbox{if }k>n.
\end{cases}
\]
It is easily seen that $|g_{n}(k)|\leq1$ (recall $\alpha\in(0,1)$ and
$z\in\lbrack0,1]$) and that
\[
g(k)=\lim_{n\rightarrow\infty}g_{n}(k)=%
\begin{cases}
1 & \mbox{if }k=0\\
(z-1)^{k} & \mbox{if }1\leq k\leq n\\
0 & \mbox{if }k>n.
\end{cases}
\]
Rewriting (\ref{temp_4}) in terms of the discrete integral on the measure
space $(\mathbb{N},\mathcal{P}(\mathbb{N}),meas_{0})$ and calling on the
Dominated Convergence Theorem, we have
\[
\varphi_{X}(z)=\lim_{n\rightarrow\infty}\int_{\mathbb{N}}g_{n}%
(k)\,meas(dk)=\int_{\mathbb{N}}g(k)\,meas(dk),
\]
which is precisely (\ref{Bern1}).

Lastly, we prove part (3). We need to show that $\varphi_{X}(z)$ admits the
representation (\ref{pcpd_pgf}). Recall $\varphi_{X}(z)=\exp\{-\phi(z)\}$ with
$\phi(z)$ of (\ref{log_Bern}). Note that%
\[
\varphi_{X}(z)=\exp\left\{  \phi(0)\left(  \frac{\phi(0)-\phi(z)}{\phi
(0)}-1\right)  \right\}
\]
and%
\begin{align*}
\frac{\phi(0)-\phi(z)}{\phi(0)}  &  =\frac{1}{\phi(0)}\sum_{n=1}^{\infty
}{\frac{p^{n}}{n(1-\alpha^{n})}}\left\{  1-(1-z)^{n}\right\} \\
&  =\frac{1}{\phi(0)}\sum_{n=1}^{\infty}\left(  -1\right)  ^{n+1}%
a_{n}(p,\alpha)z^{n},
\end{align*}
where%
\[
a_{n}(p,\alpha)=\frac{\left(  -1\right)  ^{n+1}}{\phi(0)}\sum_{j=n}^{\infty
}\binom{j}{n}{\frac{p^{j}}{j(1-\alpha^{j})}.}%
\]
Hence,%
\[
\varphi_{X}(z)=\exp\left\{  \phi(0)\left(  \sum_{n=1}^{\infty}a_{n}%
(p,\alpha)z^{n}-1\right)  \right\}  .
\]
Clearly%
\[
\sum_{n=1}^{\infty}a_{n}(p,\alpha)=1.
\]
Assume that $p<0.5$ and note that%
\begin{align*}
\sum_{n=1}^{\infty}\left\vert a_{n}(p,\alpha)\right\vert  &  \leq\frac{1}%
{\phi(0)}\sum_{n=1}^{\infty}\sum_{j=n}^{\infty}\binom{j}{n}{\frac{p^{j}%
}{j(1-\alpha^{j})}}\\
&  =\frac{1}{\phi(0)}\sum_{j=1}^{\infty}\frac{p^{j}}{j(1-\alpha^{j})}%
\sum_{n=1}^{j}\binom{j}{n}\\
&  =\frac{1}{\phi(0)}\sum_{j=1}^{\infty}\frac{p^{j}}{j(1-\alpha^{j})}\left(
2^{j}-1\right) \\
&  =\frac{1}{\phi(0)}\sum_{j=1}^{\infty}\frac{\left(  2p\right)  ^{j}%
}{j(1-\alpha^{j})}-1<\infty.
\end{align*}
Noting $\phi(0)=\lambda$ (of (\ref{pcpd_2})). $\hfill\blacksquare$

\textbf{Remark}

We note that Lemma 2 and the representation (\ref{temp_4}) of the pgf of the
Poissonian Binomial distribution (\ref{pmf_Po_Bin}) were known to Kemp (1987),
but not stated in her paper. We chose to include the technical details for the
sake of completeness. For example, the pmf of the Poissonian Binomial
distribution can also be obtained by applying the Binomial theorem to
$(z-1)^{k}$ in (\ref{temp_3}) and switching the order of summation. The
factorial moments and factorial cumulants of (\ref{pmf_Po_Bin}) can be
obtained via the expansions of $\varphi_{n-1}(1+t)$ and $\ln\varphi
_{n-1}(1+t)$ in the same way that led to (\ref{fact_mom_Bern}) and
(\ref{fact_cum_Bern}).


\begin{thebibliography}{99}                                                                                               %


\bibitem {}Abramowitz, M. and Stegun, I.A. (1965). \textit{Handbook of
Mathematical Functions}, New York: Dover.

\bibitem {}Al-Osh, M.A. and Alzaid, A.A. (1987). First-order integer valued
autoregressive INAR(1) process. Journal of Time Series Analysis, \textbf{8}, 261--275.

\bibitem {}Al-Osh, M.A., and Aly, E.-E.A.A. 1992. First order autoregressive
time-series with negative binomial and geometric marginals. Communications in
Statistics---Theory and Methods, \textbf{21}, 2483--92.

\bibitem {}Aly, E.-E.A.A and Bouzar, N. (1994). Explicit stationary
distributions for some Galton-Watson processes with immigration.
Communications in Statistics--Stochastic Models, \textbf{10}, 499--517.

\bibitem {}Aly, E.-E.A.A., and Bouzar, N. (2019). Expectation thinning
operators based on linear fractional probability generating functions. Journal
of the Indian Society for Probability and Statistics, \textbf{20}, 89--107.

\bibitem {}Benkherouf, L. and Bather, J.A. (1988). Oil exploration: Sequential
decisions in the face of uncertainty. Journal of Applied Probability,
\textbf{25}, 529--543.

\bibitem {}Bouzar, N. and Jayakumar, K. (2008). Time series with discrete
semi-stable marginals. Statistical Papers, \textbf{49}, 619--635.

\bibitem {}Feller, W. (1968). \textit{An Introduction to Probability Theory
and its Applications}, Vol. 1, John Wiley \& Sons Inc.

\bibitem {}Goldberg, K., Leighton, F.T., Newman, M., and Zuckerman, S.L.
(1976). Tables of binomial coefficients and Stirling numbers. J. of Research
of the National Bureau of Standards, \textbf{80B}, 99--171.

\bibitem {}Guerrero, M.B., Barreto-Souza, W., and Ombao, H. (2020).
Integer-valued autoregressive process with flexible marginal and innovation
distributions. arXiv:2004.08667 [stat.ME].

\bibitem {}Joe, H. (1996). Time series models with univariate margins in the
convolution-closed infinitely divisible class. Journal of Applied Probability,
\textbf{33}, 664--77.

\bibitem {}Joe, H. (2019). Likelihood inference for generalized Integer
autoregressive time series models. Econometrics, \textbf{7}, 1-13.

\bibitem {}Johnson, L.N., Kemp, A.W., and Kotz, S. (2005). \textit{Univariate
Discrete Distributions}, Third Ed., John Wiley \& Sons Inc.

\bibitem {}Jung, R.C., Ronning, G., Tremayne, A.R. (2005). Estimation in
conditional first order autoregression with discrete support. Statistical
Papers, \textbf{46}, 195--224.

\bibitem {}Kemp, A. W. (1987).A Poissonian binomial model with constrained
parameters. Naval Research Logistics, \textbf{34}, 853--858.

\bibitem {}Kemp, A.W. (1992). Heine-Euler extensions of the Poisson
distribution. Communications in Statistics--Theory and Methods, \textbf{21}, 571--588.

\bibitem {}Knopp, K. (1990). \textit{Theory and Applications of Infinite
Series}. Dover, New York.

\bibitem {}Latour, A. (1998). Existence and stochastic structure of a
non-negative integer-valued autoregressive process. Journal of Time Series
Analysis, \textbf{19}, 439--455.

\bibitem {}McKenzie, E. (1985). Some simple models for discrete variate time
series. Water Resources Bulletin, \textbf{21}, 645--650.

\bibitem {}McKenzie, E. (1988). Some ARMA models for dependent sequences of
Poisson counts. Advances in Applied Probability, \textbf{20}, 822--835.

\bibitem {}McKenzie, E. (2003). Discrete variate time series In:
\textit{Handbook of Statistics}; \textit{Stochastic Processes: Modelling and
Simulation}, C.R. Rao and D. Shanbhag, Eds., Elsevier Science, Amsterdam.

\bibitem {}Pedeli, X., Karlis, D. (2011). A bivariate INAR(1) process with
application. Statistical Modelling, \textbf{11}, 325--349.

\bibitem {}Ristic, M., Bakouch, H., \& Nastic, A. (2009). A new geometric
first-order integer-valued autoregressive (NGINAR(1)) process. Journal of
Statistical Planning and Inference, \textbf{139}, 2218--2226.

\bibitem {}Schweer, S. and Wei\ss , C.H. (2014). Compound Poisson INAR(1)
processes: Stochastic properties and testing for overdispersion. Computational
Statistics and Data Analysis, \textbf{77}, 267--284

\bibitem {}Schweer, S., Wichelhaus, C. (2015). Queuing systems of INAR(1)
processes with compound Poisson arrival distributions. Stochastic Models,
\textbf{31}, 618-635.

\bibitem {}Scotto, M.G., Wei{\ss}, C.H., and Gouveia, S. (2015). Thinning-based 
models in the analysis of integer-valued time series: a review. Statistical Modelling, 
\textbf{15}, 590--618.

\bibitem {}Smith, P.J. (1995). A recursive formulation of the old problem of
obtaining moment from cumulants and vice versa. American Statistician,
\textbf{49}, 217--218.

\bibitem {}Steutel, F.W. and van Harn, K. (1979). Discrete analogues of
self-decomposability and stability. Annals of Probability, \textbf{7}, 893--899.

\bibitem {}Wei{\ss }, C.H. (2013). Integer-valued autoregressive models for
counts showing underdispersion. Journal of Applied Statistics, \textbf{40}, 1931--1940.

\bibitem {}Wei\ss , C.H. (2018). \textit{An Introduction to Discrete-Valued
Time Series}. Hoboken: Wiley.

\bibitem {}Zhu, R. and Joe, H. (2003). A new type of discrete
self-decomposability and its application to continuous-time Markov processes
for modeling count data time series. Stochastic Models, \textbf{19}, 235--254.

\bibitem {}Zhu, R. and Joe, H. (2010). Negative binomial time series models
based on expectation thinning operators. Journal of Statistical Planning and
Inference, \textbf{140}, 1874--1888.
\end{thebibliography}
\end{document}